\documentclass[final]{siamltex}
\usepackage[letterpaper,text={6.5in,8.6in},centering]{geometry}

\DeclareMathAlphabet{\mathpzc}{OT1}{pzc}{m}{it}
\usepackage{epsfig,amssymb,graphicx,amsmath,mathrsfs}
\usepackage{subfigure}

\newcommand{\veps}{\varepsilon}
\newcommand{\T}{^{\mbox{\rm \tiny T}}}
\newtheorem{remark}{\bf Remark}
\newtheorem{prop}{\bf Proposition}

\newcommand{\norm}[1]{\ensuremath{\left\| #1 \right\|}}
\newcommand{\abs}[1]{\ensuremath{| #1 |}}
\newcommand{\bracket}[1]{\ensuremath{\left[ #1 \right]}}
\newcommand{\braces}[1]{\ensuremath{\left\{ #1 \right\}}}
\newcommand{\parenth}[1]{\ensuremath{\left( #1 \right)}}
\newcommand{\pair}[1]{\ensuremath{\langle #1 \rangle}}

\newcommand{\refeqn}[1]{(\ref{eqn:#1})}
\newcommand{\reffig}[1]{Fig. \ref{fig:#1}}
\newcommand{\tr}[1]{\mathrm{tr}\ensuremath{\negthickspace\bracket{#1}}}

\renewcommand{\Re}{\ensuremath{\mathbb{R}}}
\newcommand{\J}{\ensuremath{\mathbf{J}}}

\newcommand{\intp}{\ensuremath{\mathbf{i}}}

\newcommand{\hor}{\ensuremath{\mathrm{hor}}}

\newcommand{\dyn}{\ensuremath{\mathrm{dyn}}}
\newcommand{\geo}{\ensuremath{\mathrm{geo}}}

\title{Nonlinear Dynamics of the 3D Pendulum%
\thanks{NC, TL and NHM have been supported in part by NSF Grant ECS-0244977 and CMS-0555797. TL and ML have been supported in part by NSF Grant DMS-0504747 and DMS-0726263.}}

\author{Nalin A.~Chaturvedi\thanks{Department of Aerospace Engineering, The University of Michigan, Ann Arbor, Michigan 48109-2140 ({\tt nalin@umich.edu}).}%
\and Taeyoung Lee\thanks{Department of Aerospace Engineering, The University of Michigan, Ann Arbor, Michigan 48109-2140 ({\tt tylee@umich.edu}).}%
\and Melvin Leok\thanks{Department of Mathematics, Purdue University, West Lafayette, IN 47907-2067 ({\tt mleok@math.purdue.edu}).}%
\and N.~Harris McClamroch\thanks{Department of Aerospace Engineering, The University of Michigan, Ann Arbor, Michigan 48109-2140 ({\tt nhm@umich.edu}).}}

\begin{document}

\maketitle

\begin{abstract}
A 3D pendulum consists of a rigid body, supported at a fixed pivot, with three rotational degrees of freedom. The pendulum is acted on by a gravitational force. Symmetry assumptions are shown to lead to the planar 1D pendulum and to the spherical 2D pendulum models as special cases. The case where the rigid body is asymmetric and the center of mass is distinct from the pivot location leads to the 3D pendulum. Full and reduced 3D pendulum models are introduced and used to study important features of the nonlinear dynamics: conserved quantities, equilibria, invariant manifolds, local dynamics near equilibria and invariant manifolds, and the presence of chaotic motions. These results demonstrate the rich and complex dynamics of the 3D pendulum.
\end{abstract}



\pagestyle{myheadings}
\thispagestyle{plain}
\markboth{\sc Chaturvedi, Lee, Leok, McClamroch}{\sc Nonlinear Dynamics of the 3D Pendulum}

\section{Introduction}
Pendulum models have been a rich source of examples in nonlinear dynamics and in recent decades, in nonlinear control. The most common rigid pendulum model consists of a mass particle that is attached to one end of a massless, rigid link; the other end of the link is fixed to a pivot point that provides a rotational joint for the link and mass particle. If the link and mass particle are constrained to move within a fixed plane, the system is referred to as a planar 1D pendulum. If the link and mass particle are unconstrained, the system is referred to as a spherical 2D pendulum. Planar and spherical pendulum models have been studied in \cite{furuta, astrom}.

Numerous extensions of simple pendulum models have been proposed. These include various categories of elastic pendulum models and multi-body pendulum models. Interesting examples of multi-body pendulum models are: a pendulum on a cart, an acrobot, a pendubot, a pendulum actuated by a reaction wheel, the Furuta pendulum, and pendula consisting of multiple coupled bodies.

Pendulum models are useful for both pedagogical and research reasons. They represent physical mechanisms that can be viewed as simplified academic versions of mechanical systems that arise, for example, in robotics and spacecraft. In addition to their important role in illustrating the fundamental techniques of nonlinear dynamics, pendulum models have motivated new research directions and applications in nonlinear dynamics.

This paper considers a new 3D pendulum model, first introduced in \cite{cdc04}, and analyzes its nonlinear dynamical properties. This model consists of a rigid body, supported at a fixed pivot point that has three rotational degrees of freedom; it is acted on by a uniform gravity force.  Control and disturbance forces and moments are ignored in this development.

This paper arose out of our continuing research on a laboratory facility, referred to as the Triaxial Attitude Control Testbed (TACT). The TACT has been constructed to provide a testbed for a variety of physical experiments on attitude dynamics and attitude control. The most important feature of the TACT design is that it is supported by a three-dimensional air bearing that serves as an ideal frictionless pivot, allowing nearly unrestricted three degrees of rotation. The TACT has been described in several prior conference publications \cite{bmb, csmb}. Issues of nonlinear dynamics for the TACT have been treated in \cite{csmb, csm}. The present paper is partly motivated by the realization that the TACT is, in fact, a physical implementation of a 3D pendulum.

\section{Description of the 3D Pendulum}
A rigid 3D pendulum is a rigid body supported by a fixed, frictionless pivot, acted on by gravitational forces.
The supporting pivot allows the pendulum three rotational degrees of freedom. Uniform, constant gravity is assumed. The terminology 3D pendulum refers to the fact that the pendulum is a rigid body with three spatial dimensions and the pendulum has three rotational degrees of freedom.

Two reference frames are introduced. An inertial reference frame has its origin at the pivot; the first two axes lie in the horizontal plane and the third axis is vertical in the direction of gravity. A reference frame fixed to the pendulum body is also introduced.  The origin of this body-fixed frame is also located at the pivot.  In the body fixed frame, the moment of inertia of the pendulum is constant. This moment of inertia can be computed using the parallel axis theorem from the traditional moment of inertia with respect to a translated frame whose origin is located at the center of mass of the pendulum.  Since the origin of the body fixed frame is located at the pivot, principal axes with respect to this frame can be defined for which the moment of inertia is a diagonal matrix.  Note that the center of mass of the 3D pendulum may or may not lie on one of the principal axes defined in this way.

Rotation matrices can be used to describe the attitude of the rigid 3D pendulum. A rotation matrix maps a representation of vectors expressed in the body-fixed frame to a representation expressed in the inertial frame. Rotation matrices provide global representations of the attitude of the pendulum, which is why they are utilized here. Other attitude representations, such as exponential coordinates, quaternions, or Euler angles, can also be used following standard descriptions, but each of the representations has a disadvantage of introducing an ambiguity or singularity. In this paper, the configuration of the rigid pendulum is a rotation matrix $R$ in the special orthogonal group $SO(3)$ defined as
\[
SO(3) \triangleq \{R\in \mathbb{R}^{3\times 3}: RR\T = I_{3\times 3}, \mathrm{det}(R) = 1\}.
\]
The associated angular velocity, expressed in the body-fixed frame, is denoted by $\omega\in\mathbb{R}^3$. The constant inertia matrix, in the body-fixed frame, is denoted by the symbol $J$. The constant body-fixed vector from the pivot to the center of mass of the pendulum is denoted by $\rho$. The symbol $g$ denotes the constant
acceleration due to gravity.

Three categories of 3D pendulum models are subsequently introduced and studied. The ``full'' dynamics of the 3D pendulum are based on Euler's equations that include the gravity moment and the rotational kinematics, expressed in terms of the angular velocity and a rotation matrix; this model describes the dynamics that evolves on $TSO(3)$.   Since the gravity moment depends solely on the direction of gravity in the pendulum fixed frame, it is possible to obtain a reduced model expressed in terms of the angular velocity and a unit vector that defines the direction of gravity in the pendulum fixed frame; this model describes the dynamics that evolve on $TSO(3)/S^1$, and corresponds to the case of Lagrange--Poincar\'e reduction \cite{CeMaRa2001}.   Since there is a symmetry action given by a rotation about the gravity direction, Lagrange--Routh reduction \cite{marsden_reduc} leads to a reduced model that is restricted to a constant momentum surface and is expressed in terms of the unit vector that defines the direction of gravity in the pendulum fixed frame and its derivative; this model describes the dynamics that evolve on $TS^2$. Each of these 3D pendulum models provides special insight into the nonlinear dynamics. We develop each of these models in this paper, and we investigate the features of the nonlinear dynamics, namely invariants, equilibria, and stability for each model. 

\section{3D Pendulum Dynamics on $TSO(3)$}\label{sec:TSO}
The dynamics of the 3D pendulum are given by the Euler equation that includes the moment due to gravity:
\begin{equation}
J \dot\omega =J\omega \times \omega + mg \rho \times R\T e_3. \label{eq:Jw_dot}
\end{equation}
The rotational kinematics equations are
\begin{equation}
\dot{R} = R \widehat{\omega}.\label{eq:R_dot}
\end{equation}

Equations (\ref{eq:Jw_dot}) and (\ref{eq:R_dot}) define the full dynamics of a rigid pendulum on the tangent bundle $TSO(3)$. In the above equations, $e_3 = (0,0,1)\T$ is the direction of gravity in the inertial frame, so that $R\T e_3$ is the direction of gravity in the pendulum-fixed frame. The cross product notation $a \times b$ for vectors $a, b\in\mathbb{R}^3$ is
\begin{equation}
a \times b = [a_2b_3 - a_3b_2, \ a_3b_1 - a_1b_3, \ a_1b_2 - a_2b_1] =
\widehat{a}b,
\end{equation}
where, the skew-symmetric matrix $\widehat{a}$ is defined as
\begin{equation}
\widehat{a} = \left[\begin{array}{ccc} 0 & -a_3 & a_2 \\ a_3 &0& -a_1 \\ -a_2
&a_1& 0
\end{array}\right].
\end{equation}

A special case occurs if the center of mass of the rigid pendulum is located at the pivot. In this case $\rho = 0$, so that (\ref{eq:Jw_dot}) is given by Euler's equations with no gravity terms included. In the context of the rigid 3D pendulum, this is referred to as the balanced case. Since there is a large literature on Euler's equations with no gravity moment and the associated rotational kinematics, this case is not considered further in this paper. Rather, the focus of this paper is on the unbalanced case, where $\rho \neq 0$.

\subsection{Invariants of the 3D Pendulum}\label{ssec:TSOinvar}

There are two conserved quantities for the rigid 3D pendulum. First, the total energy, which is the sum of the rotational kinetic energy and the gravitational potential energy, is conserved.  In addition, there is a symmetry corresponding to rotations about the gravity direction through the pivot. This symmetry leads to conservation of the component of angular momentum about the gravity direction. These two results are summarized as follows.

\begin{prop}
The total energy
\[
E =\dfrac{1}{2}\, \omega\T J \omega - mg \rho\T R\T e_3,
\]
and the component of the angular momentum vector about the vertical axis through the pivot
\[
h = \omega\T J R\T e_3.
\]
are each constant along motions of the rigid 3D pendulum.
\end{prop}
\begin{proof}
The proof follows by showing that the time derivative of the total energy and the time derivative of the angular momentum component about the vertical axis are each identically zero.   We use (\ref{eq:Jw_dot}) and (\ref{eq:R_dot}) to compute the derivatives, yielding
\begin{equation*}
\dot E = \omega\T J \dot \omega - mg\rho\T \dot R\T e_3 = mg\,\omega\T( \rho
\times R\T e_3) + mg\rho\T(\omega \times R\T e_3) = 0,
\end{equation*}
and similarly,
\begin{align*}
\dot h &= \dot \omega\T J R\T e_3 + \omega\T J \dot R\T e_3, \\
&= (J\omega \times \omega)\T (R\T e_3) + (mg\rho \times R\T e_3)\T(R\T e_3) - \omega\T J(\omega \times R\T e_3), \\
&= (J\omega)\T (\omega \times R\T e_3) - \omega\T J(\omega \times R\T e_3) = 0.
\end{align*}
\end{proof}

Conservation of the angular momentum component about the vertical axis will be revisited in Section \ref{sec:TS} in the context of Noether's theorem, which states that the momentum map associated with the rotational symmetry about the gravity direction is conserved.

\subsection{Equilibria of the 3D Pendulum}
To further understand the dynamics of the 3D pendulum, we study its equilibria. Equating the RHS of (\ref{eq:Jw_dot}) and (\ref{eq:R_dot}) to zero yields
\begin{align}
J\omega_{\rm e} \times \omega_{\rm e} &+ mg\rho \times R_{\rm e}\T e_3 = 0, \label{eq:eq_omega}\\
&R_{\rm e}\widehat\omega_{\rm e} = 0. \label{eq:eq_R}
\end{align}
Since $R_{\rm e} \in SO(3)$ is non-singular, and $\widehat\cdot : \mathbb{R}^3 \rightarrow \mathbb{R}^{3\times 3}$ is a linear injection, $R_{\rm e}\widehat\omega_{\rm e} = 0$ if and only if $\omega_{\rm e} = 0$. Substituting $\omega_{\rm e} = 0$ in (\ref{eq:eq_omega}), we obtain
\begin{equation}
\rho \times R_{\rm e}\T e_3 = 0.\label{eq:equiR}
\end{equation}
Hence,
\begin{equation}
R_{\rm e}\T e_3 = \dfrac{\rho}{\|\rho\|},\label{eq:hng}
\end{equation}
or
\begin{equation}
R_{\rm e}\T e_3 = -\dfrac{\rho}{\|\rho\|}.\label{eq:inv}
\end{equation}

An attitude $R_{\rm e}$ is an equilibrium attitude if and only if the direction of gravity resolved in the body-fixed frame, $R_{\rm e}\T e_3$, is collinear with the vector $\rho$. If $R_{\rm e}\T e_3$ is in the same direction as the vector $\rho$, then $(R_{\rm e},0)$ is a \textit{hanging equilibrium} of the 3D pendulum; if $R_{\rm e}\T e_3$ is in the opposite direction as the vector $\rho$, then $(R_{\rm e},0)$ is an \textit{inverted equilibrium} of the 3D pendulum.

Thus, if $R_{\rm e}$ defines an equilibrium attitude for the 3D pendulum, then a rotation of the 3D pendulum about the gravity vector by an arbitrary angle is also an equilibrium. Consequently, in $TSO(3)$ there are two disjoint equilibrium manifolds of the 3D pendulum. The manifold corresponding to the first case in the above description is referred to as the hanging equilibrium manifold, since the center of mass is below the pivot for each attitude in the manifold. The manifold corresponding to the second case in the above description is referred to as the inverted equilibrium manifold, since the center of mass is above the pivot for each attitude in the manifold.

Following (\ref{eq:hng}) and (\ref{eq:inv}) and the discussion above, we define
\begin{align}
[R]_h &\triangleq \Big\{R \in SO(3): R\T e_3 = \dfrac{\rho}{\|\rho\|}\Big\},\label{eq:hng_att_mfld} \\
[R]_i &\triangleq \Big\{R \in SO(3): R\T e_3 =
-\dfrac{\rho}{\|\rho\|}\Big\},\label{eq:inv_att_mfld}
\end{align}
as the {\em hanging attitude manifold} and the {\em inverted attitude manifold}, respectively.

From (\ref{eq:hng}) and (\ref{eq:inv}),
\begin{equation}
\mathsf{H} \triangleq \Big\{(R,0) \in TSO(3) : R \in [R]_h\Big\},
\label{eq:hng_equi_mfld}
\end{equation}
is the manifold of hanging equilibria and
\begin{equation}
\mathtt{I} \triangleq \Big\{(R,0) \in TSO(3) : R \in [R]_i\Big\},
\label{eq:inv_equi_mfld}
\end{equation}
is the manifold of inverted equilibria, and these two equilibrium manifolds are clearly distinct.

\subsection{Local Analysis of the 3D Pendulum near an Equilibrium}\label{subsec:linear}
Consider a perturbation of the initial conditions from a hanging equilibrium $(R_{\rm e},0)$ of the 3D pendulum, using a perturbation parameter $\varepsilon$. Let $R^\veps(t)$ and $\omega^\veps(t)$ represent the perturbed solution, corresponding to initial conditions $R^\veps(0) = R_{\rm e}\exp{\veps\widehat{\delta\Theta}}$ and $\omega^\veps(0) = \veps\delta\omega$, where $\delta\Theta, \delta\omega \in \mathbb{R}^3$ are constant vectors. Note that if $\veps = 0$ then, $(R^0(0),\omega^0(0)) = (R_{\rm e},0)$ and hence
\begin{equation}
(R^0(t),\omega^0(t)) \equiv (R_{\rm e},0) \label{eq:unperturb}
\end{equation}
for all time $t \in \mathbb{R}$, which simply corresponds to the unperturbed equilibrium solution.

Consider the solution to the perturbed equations of motion for the 3D pendulum. This solution satisfies
\begin{align}
J \dot\omega^\veps &= J\omega^\veps \times \omega^\veps + mg\rho \times (R^\veps)\T e_3, \label{eq:Jw_dot_pert} \\
\dot R^\veps &= R^\veps(t)\widehat{\omega^\veps}. \label{eq:R_dot_pert}
\end{align}
Next, we differentiate both sides with respect to $\veps$ and substitute $\veps = 0$, yielding
\begin{align}
J \dot\omega_\veps^0 &= J\omega_\veps^0 \times \omega^0 + J\omega^0 \times \omega_\veps^0 + mg\rho \times (R_\veps^0)\T e_3, \label{eq:Jw_dot_pert1} \\
\dot R_\veps^0 &= R_\veps^0\widehat\omega^0 +
R^0\widehat\omega_\veps^0, \label{eq:R_dot_pert1}
\end{align}
where the subscripts denote derivatives. Substituting $R^0 = R_{\rm e}$ and $\omega^0 = 0$ from (\ref{eq:unperturb}) into (\ref{eq:Jw_dot_pert1}) and (\ref{eq:R_dot_pert1})
yields
\begin{align}
J \dot\omega_\veps^0 &=  mg\rho \times (R_\veps^0)\T e_3, \label{eq:w_lin1}\\
\dot R_\veps^0 &= R_{\rm e}\widehat\omega_\veps^0. \label{eq:R_lin1}
\end{align}
Now we define perturbation variables $\Delta\omega(t) \triangleq \omega_\veps^0(t)$ and $\widehat{\Delta\Theta}(t) \triangleq R_{\rm e}\T R_\veps^0(t)$. It can be shown that $\Delta\omega(t) = \Delta\dot\Theta(t)$. Thus, (\ref{eq:w_lin1}) and (\ref{eq:R_lin1}) can be written as
\begin{equation}
J\Delta\ddot\Theta - \dfrac{mg\widehat\rho^{\ 2}}{\|\rho\|}\Delta\Theta =
0.\label{eq:full_lin2}
\end{equation}

Note that (\ref{eq:full_lin2}) represents a mechanical system with mass matrix $J$, stiffness matrix $- \dfrac{mg\widehat\rho^{\ 2}}{\|\rho\|}$, but no damping.   Since $\widehat{\rho}^{\ 2}$ is a negative-semidefinite matrix with two negative eigenvalues and one zero eigenvalue, the stiffness matrix is positive-semidefinite with two positive eigenvalues and one zero eigenvalue. The zero eigenvalue corresponds to rotations about the vertical axis, for which gravity has no influence. To see this more explicitly, we next perform a transformation of variables.

As $\widehat\rho^{\ 2}$ is a rank 2, symmetric, negative-semidefinite matrix, it follows from \cite{bellman, bern_book} that one can simultaneously diagonalize $J$ and $\widehat\rho^{\ 2}$. Thus, there exists a non-singular matrix $M$ such that $J = MM\T$ and $-\displaystyle{\frac{mg}{\|\rho\|}\widehat\rho^{\ 2} = M\Lambda M\T}$, where $\Lambda$ is a diagonal matrix. Let $\Lambda = {\rm diag}(mgl_1,mgl_2,0)$, where $l_1$ and $l_2$ are positive. Define $x \triangleq M\T\Delta\Theta$. Then expressing $x = (x_1,x_2,x_3) \in \mathbb{R}^3$, equation (\ref{eq:full_lin2}) can be written as
\begin{align}
\ddot x_1 + mgl_1 x_1 = 0, \label{eq:x1_hang}\\
\ddot x_2 + mgl_2 x_2 = 0, \label{eq:x2_hang}\\
\ddot x_3 = 0. \label{eq:x3_hang}
\end{align}
Thus, the variable $x_3$ represents a perturbation in attitude of the 3D pendulum that corresponds to a rotation about the vertical axis.

Due to the presence of imaginary and zero eigenvalues of the linearized equations, no conclusion can be made about the stability of the hanging equilibrium or the hanging equilibrium manifold of the 3D pendulum. Indeed, the local structure of trajectories in an open neighborhood the equilibrium is that of a center manifold; there are no stable or unstable manifolds. We next show that the hanging equilibrium manifold of the 3D pendulum is Lyapunov stable.

\begin{prop} \label{thm:hang}
Consider the 3D pendulum model described by equations (\ref{eq:Jw_dot}) and (\ref{eq:R_dot}). Then, the hanging equilibrium manifold $\mathsf{H}$ given by (\ref{eq:hng_equi_mfld}) is stable in the sense of Lyapunov.
\end{prop}
\begin{proof}
Consider the following positive-semidefinite function on $TSO(3)$
\begin{equation}\label{energy}
V(R,\omega) =\frac{1}{2}~\omega^{\rm T}J \omega + mg(\|\rho\| - \rho^{\rm T}
R\T e_3).
\end{equation}
Note that $V(R,0)=0$ for all $(R,\omega) \in \mathsf{H}$ and $V(R,\omega) > 0$
elsewhere. Furthermore, the derivative along a solution of (\ref{eq:Jw_dot})
and (\ref{eq:R_dot}) is given by
\[
\begin{aligned}
  \dot{V}(R,\omega) &= \omega^{\rm T}J \dot{\omega} - mg\rho^{\rm T} \dot{R}\T e_3,\\
  &= \omega^{\rm T}(J\omega \times \omega + mg\rho \times R\T e_3) - mg\rho^{\rm T}(-\widehat\omega R\T e_3),\\
  &=mg\Big[\omega^{\rm T}(\rho \times R\T e_3) + \rho^{\rm T}(\omega \times R\T e_3)\Big]=0.
\end{aligned}
\]
Thus, $\dot V$ is negative-semidefinite on $TSO(3)$.  Also, every sublevel set of the function $V$ is compact. Therefore, the hanging equilibrium manifold $\mathsf{H}$ is Lyapunov stable. \hfill\end{proof}

Similarly, one can linearize the 3D pendulum dynamics about an equilibrium in the inverted equilibrium manifold. Expressing this linearization in terms of $(x_1,x_2,x_3) \in \mathbb{R}^3$ as in (\ref{eq:x1_hang})--(\ref{eq:x3_hang}), it can be shown that the linearization of the 3D pendulum about an inverted equilibrium can be written as
\begin{align}
\ddot x_1 - mgl_1 x_1 = 0, \label{eq:x1_inv}\\
\ddot x_2 - mgl_2 x_2 = 0, \label{eq:x2_inv}\\
\ddot x_3 = 0. \label{eq:x3_inv}
\end{align}
The linearization of the 3D pendulum about an inverted equilibrium results in a system that has two positive eigenvalues, two negative eigenvalues and two zero eigenvalues. Thus, the inverted equilibrium has a two dimensional stable manifold, a two dimensional unstable manifold and a two dimensional center manifold. It is clear that due to the presence of the two positive eigenvalues, the inverted equilibrium is unstable.

\begin{prop}
Consider the 3D pendulum model described by equations (\ref{eq:Jw_dot}) and (\ref{eq:R_dot}). Then, each equilibrium in the inverted equilibrium manifold $\mathtt{I}$ given by (\ref{eq:inv_equi_mfld}) is unstable.
\end{prop}

\section{Lagrange--Poincar\'e Reduced 3D Pendulum Dynamics on $TSO(3)/S^1$}\label{sec:TSOS}
\label{sec:red}

The equations of motion (\ref{eq:Jw_dot}) and (\ref{eq:R_dot}) for the 3D pendulum are viewed as a model for the dynamics on the tangent bundle $TSO(3)$ \cite{bloch}; these are referred to as the full equations of motion since they characterize the full attitude of the rigid pendulum. As there is a rotational symmetry corresponding to the group of rotations about the vertical axis through the pivot and the associated angular momentum component is conserved, it is possible to obtain a lower dimensional reduced model for the rigid pendulum. This Lagrange--Poincar\'e reduction is based on the fact that the dynamics and kinematics equations can be written in terms of the reduced attitude vector $\Gamma = R\T e_3 \in S^2$, which is the unit vector that expresses the gravity direction in the body-fixed frame \cite{marsden_geo_top}.

Specifically, denote the group action of $\theta\in S^1$ on $SO(3)$ by $\Phi_\theta : SO(3) \rightarrow SO(3)$, $\Phi_\theta(R) = \exp (\theta\widehat e_3) R$. This induces an equivalence class by identifying elements of $SO(3)$ that  belong to the same orbit; explicitly, for $R_1,R_2\in SO(3)$, we write $R_1\sim R_2$ if there exists a $\theta\in S^1$ such that $\Phi_\theta(R_1)=R_2$. The {\em orbit space} $SO(3)/S^1$ is the set of equivalence classes,
\begin{equation}
[R] \triangleq \{\Phi_\theta(R) \in SO(3) : \theta \in
S^1\}.\label{eq:equiv_set}
\end{equation}
For this equivalence relation, it is easy to see that $R_1 \sim R_2$ if and only if $R_1\T e_3 = R_2\T e_3$ and hence the equivalence class in (\ref{eq:equiv_set}) can alternately, be expressed as
\begin{equation}
[R] \triangleq \{R_s \in SO(3) : R_s\T e_3 = R\T e_3  \}.\label{eq:equiv}
\end{equation}
Thus, for each $R \in SO(3)$, $[R]$ can be identified with $\Gamma = R\T e_3 \in S^2$ and hence $SO(3)/S^1 \cong S^2$. This group action induces a projection $\Pi : SO(3) \rightarrow SO(3)/S^1\cong S^2$ given by $\Pi(R)=R^T e_3$.

\begin{prop}[\cite{hang_main}]\label{thm:3D_Pend_Reduced}
The dynamics of the 3D pendulum given by (\ref{eq:Jw_dot}) and (\ref{eq:R_dot})
induce a flow on the quotient space $TSO(3)/S^1$, through the projection $\pi : TSO(3) \rightarrow TSO(3)/S^1$ defined as $\pi(R,\Omega)=(R^T e_3,\Omega)$, given by the dynamics
\begin{equation}
J \dot\omega =J\omega \times \omega + mg \rho \times \Gamma,
\label{eq:Reduc_Jw_dot}
\end{equation}
and the kinematics for the reduced attitude
\begin{equation}
\dot\Gamma = \Gamma \times \omega. \label{eq:Gamma_dot}
\end{equation}
Furthermore, $TSO(3)/S^1 \cong S^2 \times \mathbb{R}^3$.
\end{prop}

Equations (\ref{eq:Reduc_Jw_dot}) and (\ref{eq:Gamma_dot}) are expressed in a
non-canonical form; they are referred to as the reduced attitude dynamics of
the 3D pendulum on $TSO(3)/S^1$.

\subsection{Special Cases of the 3D Pendulum}
Three interesting special cases are now examined. Suppose that the 3D pendulum is axisymmetric, that is two of the principal moments of inertia of the pendulum are identical and the pivot is located on the axis of symmetry of the pendulum. Assume that the body-fixed axes are selected so that $J = {\rm diag}(J_t,J_t,J_a)$ and $\rho = \rho_s e_3$ where $\rho_s$ is a positive scalar constant. Consequently, equation (\ref{eq:Reduc_Jw_dot}) can be
written in scalar form, as
\begin{equation}\label{eq:3D_spec1}
\left \{ \begin{aligned}
J_t\,\dot{\omega}_x &= (J_t - J_a)\,\omega_y\omega_z - mg\rho_s\Gamma_y, \\
J_t\,\dot{\omega}_y &= (J_a - J_t)\,\omega_z\omega_x + mg\rho_s\Gamma_x, \\
J_a\,\dot{\omega}_z &= 0.
\end{aligned}
\right.
\end{equation}
From the last equation in (\ref{eq:3D_spec1}), we see that the component of the 3D pendulum angular velocity vector about its axis of symmetry is constant. This means that a constant $\omega_z$ defines an invariant manifold for the 3D pendulum dynamics.   The special case that $\omega_z = c$, where $c \in \mathbb{R}$, leads to invariant dynamics of the axisymmetric 3D pendulum, described as follows:

\begin{prop}
Assume the 3D pendulum has a single axis of symmetry and the pivot and the center of mass are located on the axis of symmetry of the pendulum as above. The equations of motion of the 3D pendulum define an induced flow on the fiber bundle $S^2 \times \mathbb{R}^2$ corresponding to $\omega_z = c$, where $c \neq 0$, given by the equations
\begin{align*}
J_t\,\dot{\omega}_x &= c(J_t - J_a)\,\omega_y - mg\rho_s\Gamma_y, \\
J_t\,\dot{\omega}_y &= c(J_a - J_t)\,\omega_x + mg\rho_s\Gamma_x, \\
\dot{\Gamma} &= \Gamma \times \begin{bmatrix}\omega_x & \omega_y & c
\end{bmatrix}^{\rm T}.
\end{align*}
\end{prop}
\noindent These equations describe the dynamics of a {\em spinning top} where $c \in
\mathbb{R}$ denotes the spin rate.

If the spin rate of the 3D pendulum about its axis of symmetry is zero, the following result is obtained.
\begin{prop}\label{prop:spherical_pendulum}  Assume the 3D pendulum has a single axis of symmetry and the pivot and the center of mass are located on the axis of symmetry of the pendulum as above. The equations of motion of the 3D pendulum define an induced flow on the fiber bundle $S^2 \times \mathbb{R}^2$ corresponding to $\omega_z = 0$, given by the equations
\begin{align*}
J_t\,\dot{\omega}_x &=  - mg\rho_s\Gamma_y, \\
J_t\,\dot{\omega}_y &=  mg\rho_s\Gamma_x, \\
\dot{\Gamma} &= \Gamma \times \begin{bmatrix}\omega_x & \omega_y & 0
\end{bmatrix}^{\rm T}.
\end{align*}
\end{prop}
These equations describe the dynamics of a 2D {\em spherical pendulum}.

Now assume the spin rate of the 3D pendulum about its axis of symmetry is zero and, in addition, $\omega_x(0)=0$ and $\Gamma_y(0) = 0$ for the axisymmetric 3D pendulum. Then Proposition \ref{prop:spherical_pendulum} yields invariant dynamics of the axisymmetric 3D pendulum such that $\omega_x(t) = \Gamma_y(t) = 0$ for all $t \geq 0$. Therefore $\Gamma$ can be parameterized by an angle $\theta$ as $\Gamma = [-\sin\theta\ 0\ \cos\theta]\T$. This yields the following result.

\begin{prop}
Assume the 3D pendulum has a single axis of symmetry and the pivot is located on the axis of symmetry of the pendulum as above. The equations of motion of the rigid pendulum define an induced flow on the tangent bundle $TS^1$ corresponding to $\omega_x = 0$, $\omega_z = 0$ given by the equations
\begin{eqnarray*}
J_t\,\dot{\omega}_y &=& - mg\rho_s\sin\theta, \\
\dot{\theta} &=& \omega_y.
\end{eqnarray*}
\end{prop}
\noindent These equations describe the dynamics of a 1D {\em planar pendulum}.

Thus for an axially symmetric 3D pendulum with the pivot located on the axis of symmetry, the well known 2D spherical pendulum and the 1D planar pendulum can be viewed as special cases of the 3D pendulum dynamics. For an axially symmetric 3D pendulum with the pivot located on the axis of symmetry, the induced dynamics corresponding to a nonzero constant value of $\omega_z$ is fundamentally different from the dynamics of the spherical pendulum; these dynamics seem not to have been previously studied. It should be emphasized that if the 3D pendulum is asymmetric then the dynamics are general in the sense that neither the 2D spherical pendulum dynamics nor the 1D planar pendulum dynamics are special cases.

\subsection{Invariants of the Lagrange--Poincar\'e Reduced Model}

In a previous section, we obtained two integrals of motion for the full model of the 3D pendulum. In this section we summarize similar results for the Lagrange--Poincar\'e reduced model of the 3D pendulum considered here.

\begin{prop}
The total energy
\begin{align}
E =\dfrac{1}{2}\, \omega\T J \omega - mg \rho\T \Gamma,\label{eqn:E_TSOS}
\end{align}
and the component of the angular momentum vector about the vertical axis through the pivot
\[
h = \omega\T J\, \Gamma.
\]
are each constant along trajectories of the Lagrange--Poincar\'e reduced model of the 3D pendulum given by (\ref{eq:Reduc_Jw_dot}) and (\ref{eq:Gamma_dot}).
\end{prop}

\subsection{Equilibria of the Lagrange--Poincar\'e Reduced Model}
We study the equilibria of the Lagrange--Poincar\'e reduced equations of motion of the 3D pendulum given by (\ref{eq:Reduc_Jw_dot}) and (\ref{eq:Gamma_dot}). Equating the RHS of
(\ref{eq:Reduc_Jw_dot}) and (\ref{eq:Gamma_dot}) to zero yields that a natural
equilibrium $(\Gamma_{\rm e}, \omega_{\rm e})$ satisfies
\begin{align}
J\omega_{\rm e} \times \omega_{\rm e} &+ mg\rho \times \Gamma_{\rm e} = 0, \label{eq:eq_Reduc_Jomega}\\
&\Gamma_{\rm e} \times \omega_{\rm e} = 0. \label{eq:eq_Gamma}
\end{align}
Equation (\ref{eq:eq_Gamma}) implies that $\omega_{\rm e} = k\,\Gamma_{\rm e}$,
where $k \in \mathbb{R}$, and substituting this into (\ref{eq:eq_Reduc_Jomega})
yields
\begin{equation}
k^2\,J\Gamma_{\rm e} \times \Gamma_{\rm e} + mg\rho \times \Gamma_{\rm e} = 0.
\label{eq:equi_Gammae}
\end{equation}
Note that depending on whether $k$ is equal to zero or not, one can obtain {\em static} ($\omega_{\rm e} \equiv 0$) or {\em dynamic} ($\omega_{\rm e}
\not\equiv 0$) equilibrium. The dynamic case corresponds to relative equilibria \cite{HeGaLaMa2005} of the 3D pendulum.

Without loss of generality, we assume that the moment of inertia matrix is diagonal, i.e. $J=\mathrm{diag}(J_1,J_2,J_3)$, where $J_1 \geq J_2 \geq J_3 > 0$. This is achieved by choosing the body-fixed reference frame such that the body-fixed axes lie along the principal axes. Note that the vector from the pivot point to the center of mass $\rho$ may not lie along any principal axis. The following result describes the equilibria structure of the Lagrange--Poincar\'e reduced equations.

\begin{prop}\label{prop:equi_reduced}
Consider the Lagrange-Poincar\'{e} model of the 3D pendulum given by (\ref{eq:Reduc_Jw_dot}) and (\ref{eq:Gamma_dot}). The equilibria $(\Gamma_{\rm e},\omega_{\rm e})$ of the Lagrange-Poincar\'{e} model are given as follows.
\begin{enumerate}
\item The hanging equilibrium: $\parenth{\frac{\rho}{\norm{\rho}},0},$
\item The inverted equilibrium: $\parenth{-\frac{\rho}{\norm{\rho}},0},$
\item Two relative equilibria:
\begin{align}
    \braces{\parenth{-\frac{J^{-1}\rho}{\norm{J^{-1}\rho}},%
    \sqrt{\frac{mg}{\norm{J^{-1}\rho}}}J^{-1}\rho},%
    \parenth{-\frac{J^{-1}\rho}{\norm{J^{-1}\rho}},%
    -\sqrt{\frac{mg}{\norm{J^{-1}\rho}}}J^{-1}\rho}},\label{eqn:RelEqu3}
\end{align}
\item One dimensional relative equilibrium manifolds:
\begin{align}
    \parenth{-\mathrm{sgn}(\alpha)\frac{n_\alpha}{\norm{n_\alpha}},%
    \sqrt{\frac{mg}{\norm{n_\alpha}}}n_\alpha},\label{eqn:RelEqu41}
\end{align}
where $n_\alpha= (J-\frac{1}{\alpha}I_{3 \times 3})^{-1}\rho\in\Re^3$ corresponding to $\alpha \in \mathcal{L}_i$, $i \in \{1,2,3,4\}$ and
\[
\mathcal{L}_1 = (-\infty,0)\cup (\tfrac{1}{J_3},\infty),\ \mathcal{L}_2 = (0,\tfrac{1}{J_1}),\ \mathcal{L}_3 = (\tfrac{1}{J_1},\tfrac{1}{J_2}),\ \mathcal{L}_4 = (\tfrac{1}{J_2},\tfrac{1}{J_3}).
\]
Also, $\mathrm{sgn}(\cdot)$ denotes the sign function.

\item One dimensional relative equilibrium manifolds:
\begin{align}
    \parenth{-\mathrm{sgn}(\alpha)\frac{n_\alpha}{\norm{n_\alpha}},%
    -\sqrt{\frac{mg}{\norm{n_\alpha}}}n_\alpha},\label{eqn:RelEqu42}
\end{align}
where $n_\alpha$ and $\alpha$ are as given above in Case 4.
\end{enumerate}
The families of equilibria given in (\ref{eqn:RelEqu41}) and (\ref{eqn:RelEqu42}) converge to the hanging equilibrium and the inverted equilibrium when $\alpha\rightarrow 0$, and they converge to the third equilibria given in (\ref{eqn:RelEqu3}), when $\alpha\rightarrow\pm\infty$. If the vector from the pivot to the center of mass $\rho$, lies on a principal axis, i.e. $\rho\times e_i=0$ for some $i\in\{1,2,3\}$, then (\ref{eqn:RelEqu41}) and (\ref{eqn:RelEqu42}) can be rewritten as $\braces{(e_i,\gamma e_i),(-e_i,\gamma e_i)}$ for $\gamma\in\Re$.

Furthermore, there exist additional equilibria under the following assumptions on the moment of inertia matrix $J$ and the vector from the pivot to the center of mass $\rho = [\rho_1\ \rho_2\ \rho_3]\T$.
\begin{enumerate}\setcounter{enumi}{5}
\item $J_1$, $J_2$ and $J_3$ are distinct and $\rho_i = 0$. Then there exist one-dimensional relative equilibrium manifolds:
\begin{align}
\braces{\parenth{-\dfrac{p_i}{\|p_i\|},\sqrt{\dfrac{mg}{\|p_i\|}}p_i},%
\parenth{-\dfrac{p_i}{\|p_i\|},-\sqrt{\dfrac{mg}{\|p_i\|}}p_i}}
\end{align}
for $i \in \{1,2,3\}$, where %
$p_1 = \big(\gamma,\frac{\rho_2}{J_2-J_1},\frac{\rho_3}{J_3-J_1}\big)$, %
$p_2 = \big(\frac{\rho_1}{J_1-J_2},\gamma,\frac{\rho_3}{J_3-J_2}\big)$, %
$p_3 = \big(\frac{\rho_1}{J_1-J_3},\frac{\rho_2}{J_2-J_3},\gamma\big)$ and $\gamma \in \mathbb{R}$.
\item $J_1 = J_2 \not= J_3$.
    \begin{enumerate}
    \item If $\rho_1 = \rho_2 = 0$. Then there exist two-dimensional relative equilibrium manifolds:
    \[
    \bigg\{\left(-\dfrac{q}{\|q\|},\sqrt{\dfrac{mg}{\|q\|}}q\right),%
           \left(-\dfrac{q}{\|q\|},-\sqrt{\dfrac{mg}{\|q\|}}q\right)\bigg\}
    \]
    where $q = \big(\gamma,\delta,\frac{\rho_3}{J_3-J_1}\big)$ and $\gamma,\delta \in \mathbb{R}$.
    \item If $\rho_3 = 0$. Then there exist one-dimensional relative equilibrium manifolds:
    \[
    \bigg\{\left(-\dfrac{r}{\|r\|},\sqrt{\dfrac{mg}{\|r\|}}r\right),%
    \left(-\dfrac{r}{\|r\|},-\sqrt{\dfrac{mg}{\|r\|}}r\right)\bigg\}
    \]
    where $r = \big(\frac{\rho_1}{J_1-J_3},\frac{\rho_2}{J_1-J_3},\gamma\big)$ and $\gamma \in \mathbb{R}$.
    \end{enumerate}
\end{enumerate}
\end{prop}
\begin{proof}
From (\ref{eq:equi_Gammae}), an equilibrium $(\Gamma_{\rm e},\omega_{\rm e})$ satisfies \begin{align}
    k^2 J\Gamma_{\rm e} + m g \rho  =k_1\Gamma_{\rm e},\label{eqn:cond0}
\end{align}
for a constant $k_1\in\Re$. We solve this equation to obtain the expression for an equilibrium attitude $\Gamma_{\rm e}$ for two cases; when $k_1 =0$ and when $k_1\neq 0$. The corresponding value of the constant $k$ yields the expression for the equilibrium angular velocity as $\omega_{\rm e}=k\Gamma_{\rm e}$.

\textit{Equilibria 3:} Suppose $k_1 =0$. It follows that $k\neq 0$ from \refeqn{cond0}. Thus, we have $\Gamma_{\rm e}=-\frac{mg}{k^2}J^{-1}\rho$. Since $\norm{\Gamma_{\rm e}}=1$, we obtain $k^2=mg \|J^{-1}\rho\|$, which gives \refeqn{RelEqu3}.

\textit{Equilibria 1, 2, 4 and 5:} Suppose $k_1\neq 0$. If $k = 0$, \refeqn{cond0} yields the hanging and the inverted equilibrium. Suppose $k \not= 0$. Define $\alpha=\frac{k^2}{k_1}\in\Re\backslash\{0\}$, and $v=k_1\Gamma_{\rm e}\in\Re^3$. Then, \refeqn{cond0} can be written as
\begin{align}
    (\alpha J - I_{3\times 3})v = -mg \rho.\label{eqn:cond1}
\end{align}
Note that for $\alpha\in\Re\backslash\{0,\frac{1}{J_1},\frac{1}{J_2},\frac{1}{J_3}\}$ the matrix $(J - \tfrac{1}{\alpha}I_{3\times 3})$ is invertible. Then, \refeqn{cond1} can be solved to obtain $v=-\frac{mg}{\alpha}(J - \tfrac{1}{\alpha}I_{3\times 3})^{-1} \rho$. Since $\norm{\Gamma_{\rm e}}=1$, we have $\norm{v}=\norm{k_1 \Gamma_{\rm e}}=\abs{k_1}$. We consider two sub-cases; when $k_1 >0$, and $k_1 <0$.

If $k_1 > 0$, we have $k_1=\norm{v}$ and $\alpha > 0$. Thus, we obtain the expression for  equilibria attitudes as
\begin{align}
    \Gamma_{\rm e} = \frac{v}{\norm{v}} = - \frac{n_\alpha}{\norm{n_\alpha}} = -\mathrm{sgn}(\alpha) \frac{n_\alpha}{\norm{n_\alpha}}.\label{eqn:Game4}
\end{align}
where $n_\alpha = (J - \tfrac{1}{\alpha}I_{3\times 3})^{-1} \rho\in\Re^3$. Since $k^2 = \alpha k_1=\alpha\norm{v}$, we obtain the expression for equilibria angular velocities as
\begin{align}
    \omega_{\rm e} = k\Gamma_{\rm e} =  \mp \sqrt{\alpha\norm{v}}\frac{n_\alpha}{\norm{n_\alpha}}=\mp \sqrt{\frac{mg}{\norm{n_\alpha}} }n_\alpha.\label{eqn:Omegae4}
\end{align}
Thus, \refeqn{Game4} and \refeqn{Omegae4} correspond to the families of equilibria given by \refeqn{RelEqu41} and \refeqn{RelEqu42} for $\alpha > 0$. Consider the limiting case when $\alpha \rightarrow \infty$. We have
\begin{align*}
    \lim_{\alpha\rightarrow\infty} - \frac{n_\alpha}{\norm{n_\alpha}} = \lim_{\alpha\rightarrow\infty} -\frac{(J-I_{3\times3}/\alpha)^{-1}\rho}{\norm{(J- I_{3\times3}/\alpha)^{-1}\rho}} = -\frac{J^{-1}\rho}{\norm{J^{-1}\rho}}.
\end{align*}
Similarly,
\begin{align*}
\lim_{\alpha\rightarrow\infty} \sqrt{\frac{mg}{\norm{n_\alpha}}}n_\alpha = \sqrt{\frac{mg}{\norm{J^{-1}\rho}}}J^{-1}\rho.
\end{align*}
Thus, as $\alpha\rightarrow\infty$, \refeqn{Game4} and \refeqn{Omegae4} converges to the first relative equilibria given in \refeqn{RelEqu3}.

Similarly, if $k_1 < 0$, we have $k_1 = -\norm{v}$, $\alpha < 0$, and $k^2 = -\alpha\norm{v}$. Thus, the relative equilibria are described by
\begin{gather}
\Gamma_{\rm e}= - \frac{v}{\|v\|} = \frac{n_\alpha}{\norm{n_\alpha}} = -\mathrm{sgn}(\alpha)\frac{n_\alpha}{\norm{n_\alpha}},\quad \omega_{\rm e}=\pm\sqrt{\frac{mg}{\norm{n_\alpha}}}n_\alpha,
\end{gather}
which corresponds to the families of equilibria given by \refeqn{RelEqu41} and \refeqn{RelEqu42} for $\alpha < 0$. It can be similarly shown that they converges to the third relative equilibria given by \refeqn{RelEqu3} as $\alpha\rightarrow-\infty$.

\par
\par
Next, consider \refeqn{RelEqu41} and \refeqn{RelEqu42} as $\alpha \rightarrow 0$. Expressing $n_\alpha = \alpha(\alpha J-I_{3\times 3})^{-1}\rho$, we obtain
\[
\lim_{\alpha \rightarrow 0^+} -\mathrm{sgn}(\alpha)\frac{n_\alpha}{\|n_\alpha\|} = \lim_{\alpha \rightarrow 0^+} -\frac{(\alpha J-I_{3\times 3})^{-1}\rho}{\|(\alpha J-I_{3\times 3})^{-1}\rho\|} = -\frac{\rho}{\|\rho\|},
\]
which corresponds to the inverted attitude. Similarly,
\[
\lim_{\alpha \rightarrow 0^-} -\mathrm{sgn}(\alpha)\frac{n_\alpha}{\|n_\alpha\|} = \lim_{\alpha \rightarrow 0^-} \frac{(\alpha J-I_{3\times 3})^{-1}\rho}{\|(\alpha J-I_{3\times 3})^{-1}\rho\|} = \frac{\rho}{\|\rho\|},
\]
which corresponds to the hanging attitude. For the angular velocity term,
\[
\lim_{\alpha \rightarrow 0} \sqrt{\dfrac{mg}{\|n\|}}n = \lim_{\alpha \rightarrow 0} \sqrt{mg\|n\|}\frac{n}{\|n\|} = \pm\frac{\rho}{\|\rho\|} \lim_{\alpha \rightarrow 0} \sqrt{mg\|\alpha\rho\|} = 0.
\]
Thus, as $\alpha \rightarrow 0^-$, \refeqn{RelEqu41} and \refeqn{RelEqu42} yields the hanging equilibrium given in case (1), and similarly, as $\alpha \rightarrow 0^+$, \refeqn{RelEqu41} and \refeqn{RelEqu42} yields the inverted equilibrium given in case (2).

\par
Now suppose the vector from the pivot to the center of mass $\rho$, lies on a principal axis, i.e. $\rho\times e_i=0$ for some $i\in\{1,2,3\}$. Then the vector $\rho$ can be expressed as $\rho = s e_i$, where $s \in \mathbb{R}$. Then for all $\alpha \in \Re\backslash\{0,\frac{1}{J_1},\frac{1}{J_2},\frac{1}{J_3}\}$, $(\alpha J-I_{3 \times 3})$ is an invertible diagonal matrix, and hence $n_\alpha = \frac{1}{\alpha}(\alpha J - I_{3 \times 3})^{-1}\rho = \dfrac{s}{\alpha(\alpha J_i - 1)} e_i$. Then, it follows from (\ref{eqn:RelEqu41}) and (\ref{eqn:RelEqu42}) that the equilibria can be rewritten as $\braces{(e_i,\gamma e_i),(-e_i,\gamma e_i)}$ for $\gamma\in\Re$.

The solution of (\ref{eqn:cond1}) for $\alpha \in \mathbb{R}$ yields all possible equilibria of (\ref{eq:Reduc_Jw_dot}) and (\ref{eq:Gamma_dot}). Equation (\ref{eqn:RelEqu41}) and (\ref{eqn:RelEqu42}) present the solution of (\ref{eqn:cond1}) for all $\alpha \in \Re\backslash\{0,\frac{1}{J_1},\frac{1}{J_2},\frac{1}{J_3}\}$. As shown before, $\alpha = 0$ yields the hanging and the inverted equilibrium. If $\alpha \in \{\frac{1}{J_1},\frac{1}{J_2},\frac{1}{J_3}\}$, (\ref{eqn:cond1}) can have solutions under certain conditions. This yields the additional equilibria of (\ref{eq:Reduc_Jw_dot}) and (\ref{eq:Gamma_dot}) given in Case 6 and 7.

\textit{Equilibria 6:} Suppose $J_1,J_2$ and $J_3$ are distinct. Then it is easy to see that for $\alpha = 1/J_1$, (\ref{eqn:cond1}) has a solution iff $\rho_1 = 0$. In this case, (\ref{eqn:cond1}) can be written as
\[
\begin{bmatrix} 0 & 0& 0\\0 &J_2-J_1 & 0 \\ 0 &0 &J_3-J_1 \end{bmatrix}v = -mgJ_1\begin{bmatrix} 0 \\\rho_2 \\ \rho_3 \end{bmatrix}.
\]
Since $\alpha > 0$, it can be shown as in (\ref{eqn:Game4}) that $\Gamma_{\rm e} = -\dfrac{p_1}{\|p_1\|}$ and $\omega_{\rm e} = \pm\sqrt{\dfrac{mg}{\|p_1\|}}p_1$, where $p_1 = \big(\gamma,\frac{\rho_2}{J_2-J_1},\frac{\rho_3}{J_3-J_1}\big)$ and $\gamma \in \mathbb{R}$. Similarly, one can yield solutions of (\ref{eqn:cond1}) for the case where $\alpha = 1/J_2$ and $\alpha = 1/J_3$ iff $\rho_2 = 0$ and $\rho_3 = 0$, respectively. Thus for distinct principal moments of inertia, one obtains the equilibria given in case 5, corresponding to the condition $\rho_i = 0$ where $i \in \{1,2,3\}$.

\textit{Equilibria 7:} (a) Suppose $J_1 = J_2 \not= J_3$. Then it is easy to see that for $\alpha = 1/J_1 = 1/J_2$, (\ref{eqn:cond1}) has a solution iff $\rho_1 = \rho_2 = 0$. In this case, (\ref{eqn:cond1}) can be written as
\[
\begin{bmatrix} 0 & 0& 0\\0 &0 & 0 \\ 0 &0 &J_3-J_1 \end{bmatrix}v = -mgJ_1\begin{bmatrix} 0 \\0 \\ \rho_3 \end{bmatrix}.
\]
Since $\alpha > 0$, it can be shown as in (\ref{eqn:Game4}) that $\Gamma_{\rm e} = -\dfrac{q}{\|q\|}$ and $\omega_{\rm e} = \pm\sqrt{\dfrac{mg}{\|q\|}}q$, where $q = \big(\gamma,\delta,\frac{\rho_3}{J_3-J_1}\big)$ and $\gamma,\delta \in \mathbb{R}$.

(b) Similar to above, one can yield solutions of (\ref{eqn:cond1}) for the case where $\alpha = 1/J_3$ iff $\rho_3 = 0$.

Thus for $J_1 = J_2 \not= J_3$, one obtains the equilibria given in case 7(a) and 7(b) under the specified conditions. Finally for the case where $J_1 = J_2 = J_3$, there are no additional solutions of (\ref{eqn:cond1}) for $\alpha \in \{\frac{1}{J_1},\frac{1}{J_2},\frac{1}{J_3}\}$.\hfill\end{proof}

\begin{figure}
\centerline{
\subfigure[A hanging equilibrium]%
    {\includegraphics[width=0.28\textwidth]{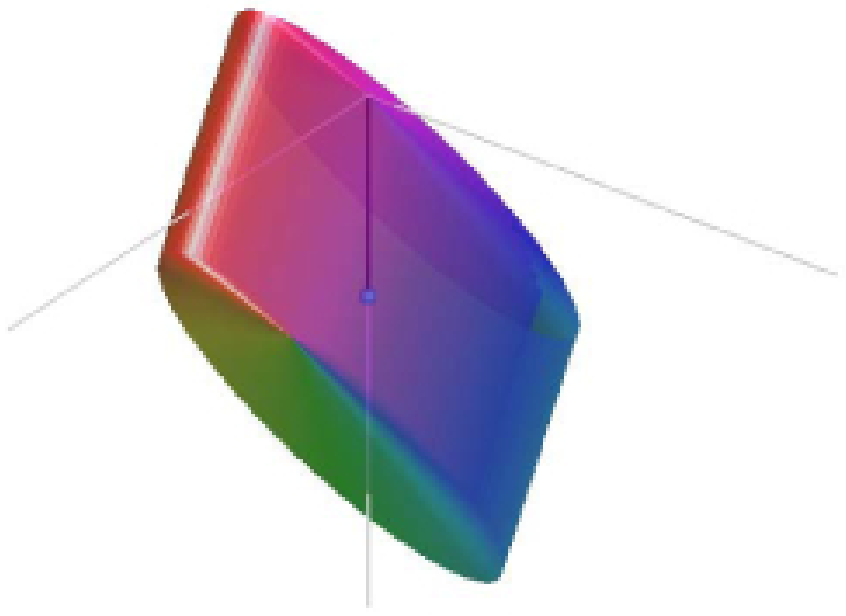}\label{fig:RelEquSphHe}}
\subfigure[3D view]%
    {\includegraphics[width=0.24\textwidth]{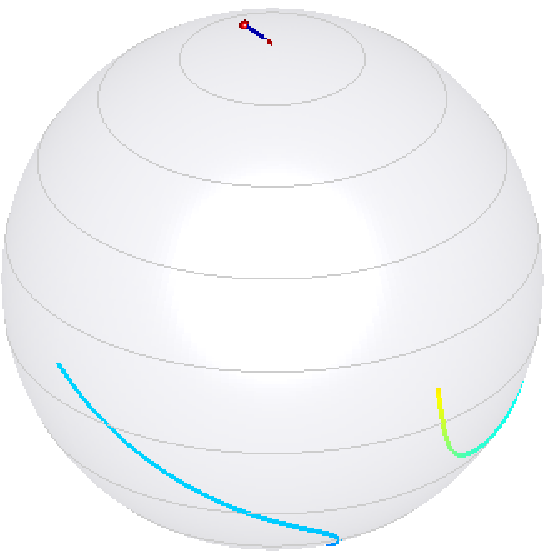}\label{fig:RelEquSph3d}}
\subfigure[Top view]%
    {\includegraphics[width=0.24\textwidth]{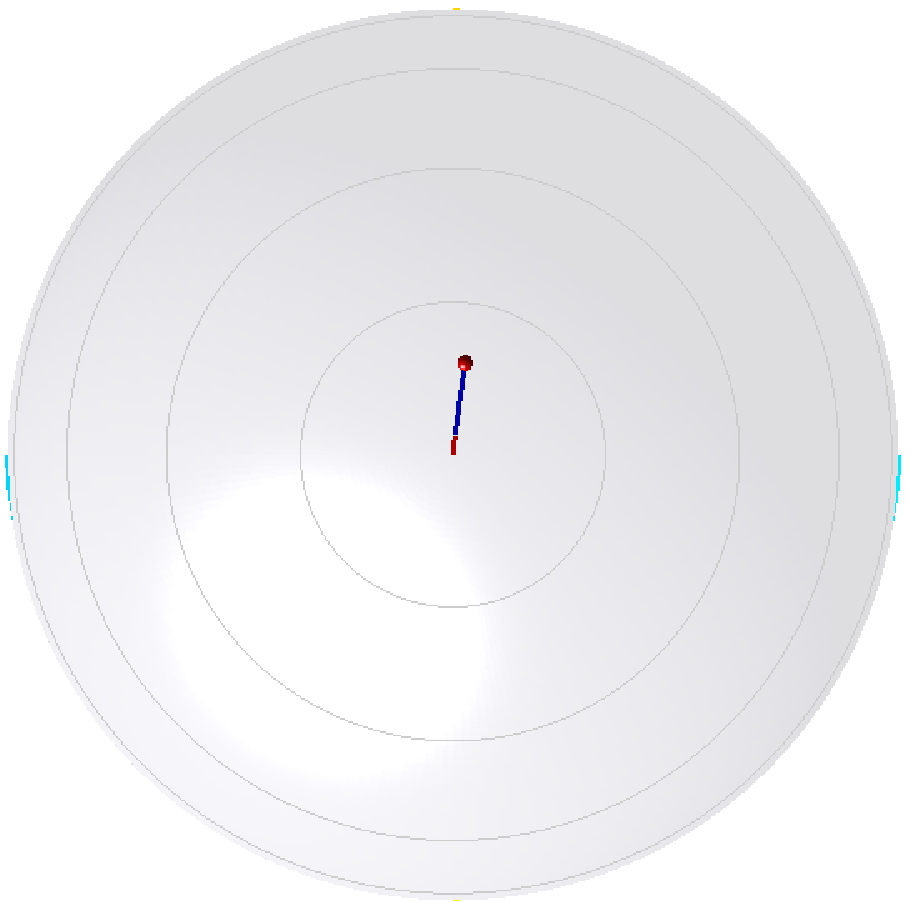}\label{fig:RelEquSphTop}}
\subfigure[Bottom view]%
    {\includegraphics[width=0.24\textwidth]{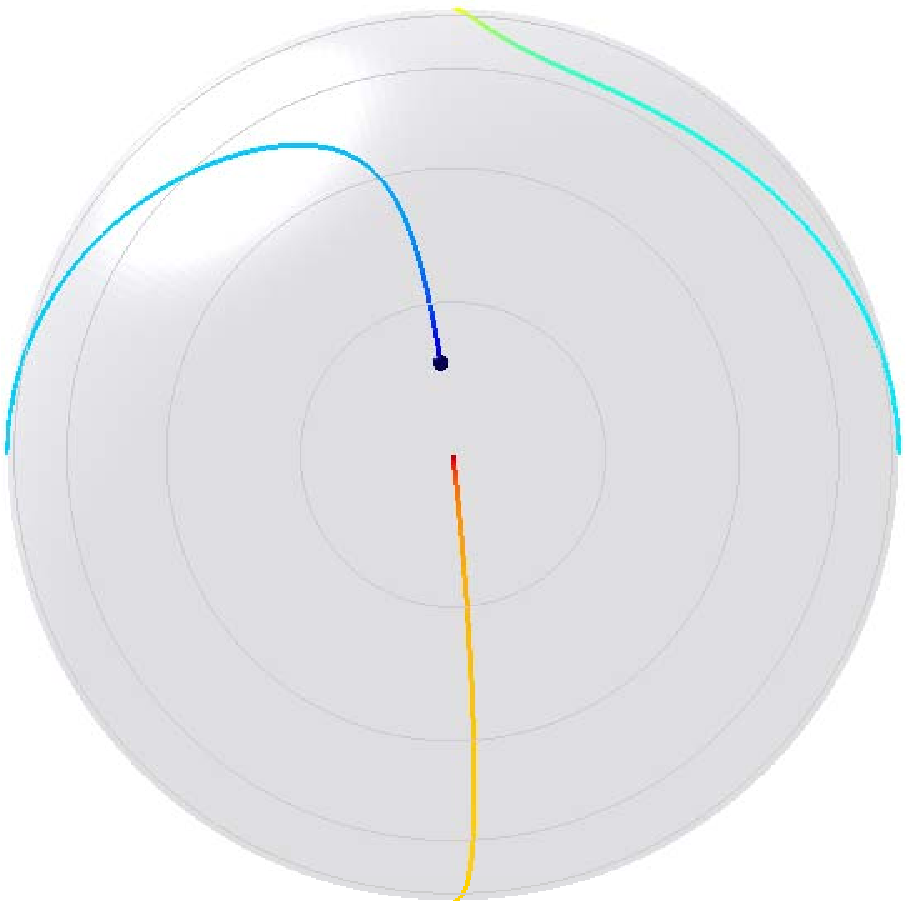}\label{fig:RelEquSphBottom}}}
\caption{Relative equilibria attitudes for an elliptic cylinder 3D pendulum model}
\end{figure}

\paragraph{Numerical example} We show the equilibrium structure of a particular 3D pendulum model. We choose an elliptic cylinder with its semimajor axis $a=0.8\,m$, semi-minor axis $b=0.2\,m$, and height $0.6\,m$. The pivot point is located at the surface of the upper ellipse, and it is offset from the center by
$[-\frac{a}{2}, \frac{b}{2}, 0]$. The moment of inertia is given by $J=\mathrm{diag}(0.4486, 0.3943,0.0772)$ and the vector from the pivot to the mass center is $\rho=[-0.0140,0.1044,0.4989]$. One of hanging attitudes is shown in \reffig{RelEquSphHe}.

Figures \ref{fig:RelEquSph3d}--\ref{fig:RelEquSphBottom} show the relative equilibria attitudes on $S^2$, where the top corresponds  $\Gamma=-e_3$, and the bottom corresponds to $\Gamma=e_3$. The inverted equilibrium is denoted by a red dot, and the hanging equilibrium is denoted by a blue dot, and the equilibria of \refeqn{RelEqu3} are located at the intersection of the blue line and the red line in \reffig{RelEquSphTop}. The families of the relative equilibria given by \refeqn{RelEqu41} and \refeqn{RelEqu42} are shown by the four segments of solid lines corresponding to $\mathcal{L}_i$, where $i \in \{1,2,3,4\}$ and the value of $\alpha$ is represented by color-shading; $\alpha$ varies from $-\infty$ (blue color) to $\infty$ (red color). Note that the reduced attitude in both \refeqn{RelEqu41} and \refeqn{RelEqu42}, are the same and these families of equilibria only differ by a sign in the angular velocity vector at the equilibrium.

The relative equilibria for $\alpha\in(-\infty,0)$ are shown by a segment of a blue line in \reffig{RelEquSphTop}, which starts from the third equilibria  given by \refeqn{RelEqu3}, and converges to the inverted equilibrium. For $\alpha\in(0,\frac{1}{J_3})$, three disjoint segments of the relative equilibria attitudes  are shown in \reffig{RelEquSphBottom}; $\alpha\in(0,\frac{1}{J_1}) = \mathcal{L}_2$ (upper-left, moving counter-clockwise from blue to cyan), $\alpha\in(\frac{1}{J_1},\frac{1}{J_2}) = \mathcal{L}_3$ (upper-right, moving counter-clockwise from cyan to green), and $\alpha\in(\frac{1}{J_2},\frac{1}{J_3}) = \mathcal{L}_4$ (lower center, moving upward from green to orange). The relative equilibria attitudes for $\alpha\in(\frac{1}{J_3},\infty)$ are shown in \reffig{RelEquSphTop} by a red line segment, which converges to the blue line at the third equilibrium given by \refeqn{RelEqu3}.

Since no component of the center of mass vector vanishes, there are no additional
equilibria.   In summary, the hanging attitude, the inverted attitude, and the attitude given by \refeqn{RelEqu3} are equilibrium attitudes, and there are four mutually disjoint equilibrium attitude segments corresponding to $\mathcal{L}_i$, where $i \in \{1,2,3,4\}$. 

\paragraph{Relation to the equilibrium manifolds of the full model} Let
\[
\Gamma_h \triangleq \dfrac{\rho}{\|\rho\|},\ {\rm and}\ \Gamma_i \triangleq
-\dfrac{\rho}{\|\rho\|}.
\]
Then it follows from Proposition \ref{prop:equi_reduced} that $(\Gamma_h,0)$ and $(\Gamma_i,0)$ are equilibria of the Lagrange--Poincar\'e reduced model of the 3D pendulum. These are called the hanging equilibrium and the inverted equilibrium of the Lagrange--Poincar\'e reduced model, respectively.

Let $(R_{\rm e},0)$ denote an equilibrium in either the hanging equilibrium manifold or the inverted equilibrium manifold of the full equations (\ref{eq:Jw_dot}) and (\ref{eq:R_dot}) and $\pi :  TSO(3) \rightarrow TSO(3)/S^1$ be the projection as in Proposition \ref{thm:3D_Pend_Reduced}. Then, it can be shown that either $\pi(R_{\rm e},0) = (\Gamma_h,0)$ or $\pi(R_{\rm e},0) = (\Gamma_i,0)$. Thus, the hanging and the inverted equilibrium manifold of the full equations are identified with the hanging and the inverted equilibrium of the Lagrange--Poincar\'e reduced equations.

To study the properties of the equilibrium manifolds, it is advantageous to consider the hanging and the inverted equilibrium of the equations of the 3D pendulum in terms of the reduced attitude as in (\ref{eq:Reduc_Jw_dot}) and (\ref{eq:Gamma_dot}). The following result provides the identification.

\begin{prop}[\cite{hang_main}]\label{prop:identif}
The hanging equilibrium manifold and the inverted equilibrium manifold of the
3D pendulum given by (\ref{eq:Jw_dot}) and (\ref{eq:R_dot}) are identified with
the hanging equilibrium $(\Gamma_h,0)$ and the inverted equilibrium
$(\Gamma_i,0)$ of the reduced attitude equations given by
(\ref{eq:Reduc_Jw_dot}) and (\ref{eq:Gamma_dot}).
\end{prop}

\subsection{Local Analysis of the Lagrange--Poincar\'e Reduced Model near an Equilibrium}

In the last section, we showed that the Lagrange--Poincar\'e reduced model of the 3D
pendulum can have equilibria with non-zero angular velocities. Also, there are only two {\em static} equilibria, namely the hanging equilibrium and the inverted equilibrium. As stated in Proposition \ref{prop:identif}, these equilibria correspond to the disjoint equilibrium manifolds of the full equations of the 3D pendulum.

We next focus on these static equilibria of the Lagrange--Poincar\'e reduced equations. The
identification mentioned in Proposition \ref{prop:identif} relates properties of the equilibrium
manifolds of the full equations and the equilibria of the Lagrange--Poincar\'e reduced equations. We deduce the stability of the hanging and the inverted equilibrium manifolds of the full equations by studying the stability property of the hanging equilibrium and the inverted equilibrium of the Lagrange--Poincar\'e reduced equations.

Consider the linearization of (\ref{eq:Reduc_Jw_dot})--(\ref{eq:Gamma_dot})
about an equilibrium $(\Gamma_h,0)= (R_{\rm e}\T e_3,0)$, where $(R_{\rm
e},0)$ is an equilibrium of the hanging equilibrium manifold $\mathsf{H}$. Since
${\rm dim}\Big[TSO(3)/S^1\Big] = 5$, the linearization of
(\ref{eq:Reduc_Jw_dot})--(\ref{eq:Gamma_dot}) about $(\Gamma_h,0)$ evolves on
$\mathbb{R}^5$.

\begin{prop}\label{prop:Reduc_3D_Pend_linear_hng}
The linearization of the Lagrange--Poincar\'e reduced equations for the 3D pendulum, about
the equilibrium $(\Gamma_h,0) = (R_{\rm e}\T e_3,0)$ described by equations
(\ref{eq:Reduc_Jw_dot})--(\ref{eq:Gamma_dot}) can be expressed using
$(x_1,x_2,\dot x_1, \dot x_2, \dot x_3) \in \mathbb{R}^5$ according to
(\ref{eq:x1_hang})--(\ref{eq:x3_hang}).
\end{prop}

\begin{proof}
Consider a perturbation in terms of
the perturbation parameter $\veps \in \mathbb{R}$ as before. Let
$(R^\veps(t),\omega^\veps(t))$ denote the perturbed solution of
(\ref{eq:Jw_dot})--(\ref{eq:R_dot}). Since $\Gamma = R^{\rm T}
e_3$, the perturbed solution of (\ref{eq:Reduc_Jw_dot})--(\ref{eq:Gamma_dot})
is given by $(\Gamma^\veps(t),\omega^\veps(t))$ where $\Gamma^\veps(t) =
(R^\veps)\T(t) e_3$. Define the perturbation variables  $\Delta\omega(t) \triangleq
\omega_\veps^0(t)$ and $\Delta\Gamma(t) \triangleq \Gamma_\veps^0(t) =
(R_\veps^0)\T(t) e_3$. From definition of $\Delta\Theta$ in Subsection
\ref{subsec:linear}, note that
\[
\Delta\Gamma = -\widehat{\Delta\Theta} R_{\rm e}\T e_3 = \widehat\Gamma_h
\Delta\Theta \in T_{\Gamma_h}S^2.
\] Then from (\ref{eq:w_lin1}) and the definition of $\Delta\Gamma$, it can be easily shown that the linearization of (\ref{eq:Reduc_Jw_dot})--(\ref{eq:Gamma_dot}) is given by
\begin{align}
J\Delta\dot\omega &=  mg\widehat\rho\ \Delta\Gamma, \label{eq:linearr1} \\
\Delta\dot\Gamma &= \widehat\Gamma_h\,\Delta\omega. \label{eq:linearr2}
\end{align}
Next, we express (\ref{eq:linearr1}) and (\ref{eq:linearr2}) in terms of
$(x,\dot x)$. Specifically, we show that $(\Delta\Gamma,\Delta\omega) \in
T_{\Gamma_h}S^2 \times \mathbb{R}^3$ can be expressed using $(x_1,x_2,\dot x_1,
\dot x_2, \dot x_3) \in \mathbb{R}^5$.

Since $x = M\T \Delta\Theta$, and $M$ is nonsingular, $\Delta\omega =
M^{-\mbox{\rm \tiny T}}\dot x$ and $\displaystyle{\Delta\Gamma =
\widehat\Gamma_h M^{-\mbox{\rm \tiny T}} x}$. We now give an orthogonal
decomposition of the vector $\Delta\Theta = M^{-\mbox{\rm \tiny T}}x$ into a
component along the vector $\rho$ and a component normal to the vector $\rho$.
This decomposition is
\[
M^{-\mbox{\rm \tiny T}} x = -\frac{\widehat\rho^{\
2}}{\|\rho\|^2}(M^{-\mbox{\rm \tiny T}} x) +
\frac{1}{\|\rho\|^2}\Big[\rho\T(M^{-\mbox{\rm \tiny T}} x)\Big]\rho,
\]
where $\dfrac{1}{\,\|\rho\|^2}\Big[\rho\T(M^{-\mbox{\rm \tiny T}} x)\Big]\rho
\in {\rm span}\{\rho\}$ and $-\dfrac{\widehat\rho^{\
2}}{\,\|\rho\|^2}(M^{-\mbox{\rm \tiny T}} x) \in {\rm span}\{\rho\}^\perp$.
Thus,
\[
\Delta\Gamma = \widehat\Gamma_h\,\Delta\Theta = \frac{\widehat\rho}{\|\rho\|}
M^{-\mbox{\rm \tiny T}} x =\frac{1}{mg\|\rho\|^2}\widehat\rho\, M\Lambda\, x,
\]
does not depend on $x_3$ since $\Lambda={\rm diag}(mgl_1,mgl_2,0)$. Thus, we
can express the linearization of (\ref{eq:Reduc_Jw_dot})--(\ref{eq:Gamma_dot})
at $(\Gamma_h,0) = (R_{\rm e}\T e_3,0)$ in terms of the variables
$(x_1,x_2,\dot x_1, \dot x_2, \dot x_3)$ according to
(\ref{eq:x1_hang})--(\ref{eq:x3_hang}).\hfill\end{proof}

\begin{remark} {\rm Note that due to our careful choice of variables, one can discard $x_3$ from (\ref{eq:x1_hang})--(\ref{eq:x3_hang}) when studying the stability properties of the inverted equilibrium manifold. Thus, $x_3$ corresponds to a component of the perturbation in the attitude that is tangential to the inverted equilibrium manifold. However, the angular velocity corresponding to $x_3$ given by $\dot{x}_3$ is retained.}
\end{remark}

Summarizing the above, the linearization of (\ref{eq:Reduc_Jw_dot})--(\ref{eq:Gamma_dot}) about the hanging equilibrium
$(\Gamma_h,0)$ is expressed as
\begin{align}
\ddot x_1 + mgl_1 x_1 = 0, \label{eq:x1_hang_red}\\
\ddot x_2 + mgl_2 x_2 = 0, \label{eq:x2_hang_red}\\
\dot x_3 = 0. \label{eq:x3_hang_red}
\end{align}
It is clear that due to the presence of zero and imaginary eigenvalues, one cannot arrive at a conclusion about the stability of the hanging equilibrium $(\Gamma_h,0)$ from the linear analysis. Therefore we next consider Lyapunov analysis.

\begin{prop}\label{prop:hang_reduc}  The hanging equilibrium $(\Gamma_h,0) = \left(\frac{\rho}{\|\rho\|},0\right)$, of the reduced dynamics of the 3D pendulum described by equations (\ref{eq:Reduc_Jw_dot}) and (\ref{eq:Gamma_dot}) is stable in the sense of Lyapunov.
\end{prop}
\begin{proof}
Consider the candidate Lyapunov function
\begin{equation}\label{energy_reduc}
  V(\Gamma,\omega) =\frac{1}{2}~\omega^{\rm T}J \omega + mg(\|\rho\| - \rho^{\rm T} \Gamma).
\end{equation}
Note that $V(\Gamma_h,0)=0$ and $V(\Gamma, \omega) > 0$ elsewhere. Furthermore,
the derivative along a solution of (\ref{eq:Reduc_Jw_dot}) and
(\ref{eq:Gamma_dot}) is given by
\[
\begin{aligned}
  \dot{V}( \Gamma, \omega) &= \omega^{\rm T}J \dot{\omega} - mg\rho^{\rm T} \dot{\Gamma},\\
  &= \omega^{\rm T}(J\omega \times \omega + mg\rho \times \Gamma) - mg\rho^{\rm T}(\Gamma \times \omega),\\
  &=\omega^{\rm T}mg\rho \times \Gamma-mg\rho^{\rm T}\Gamma \times \omega=0.
\end{aligned}
\]
Thus, the hanging equilibrium is Lyapunov stable.\hfill\end{proof}

\begin{remark}{\rm Note that combining Proposition \ref{prop:hang_reduc} with Proposition \ref{prop:identif} immediately confirms the stability result for the hanging equilibrium manifold in Proposition \ref{thm:hang}.}
\end{remark}

We next study the local properties of the Lagrange--Poincar\'e reduced equations of the 3D pendulum near the inverted equilibrium $(\Gamma_i,0)$. Consider the linearization of
(\ref{eq:Reduc_Jw_dot})--(\ref{eq:Gamma_dot}) about an equilibrium $(\Gamma_i,0)= (R_{\rm e}\T e_3,0)$, where $(R_{\rm e},0)$ is an equilibrium of the inverted equilibrium manifold $\mathtt{I}$. A result similar to Proposition \ref{prop:Reduc_3D_Pend_linear_hng} follows.

\begin{prop}\label{prop:Reduc_3D_Pend_linear_inv}
The linearization of the Lagrange--Poincar\'e reduced equations for the 3D pendulum, about the equilibrium $(\Gamma_i,0) = (R_{\rm e}\T e_3,0)$ described by equations (\ref{eq:Reduc_Jw_dot})--(\ref{eq:Gamma_dot}) can be expressed using $(x_1,x_2,\dot x_1, \dot x_2, \dot x_3) \in \mathbb{R}^5$ according to (\ref{eq:x1_inv})--(\ref{eq:x3_inv}).
\end{prop}

Summarizing the above, the linearization of (\ref{eq:Reduc_Jw_dot})--(\ref{eq:Gamma_dot}) about the inverted equilibrium $(\Gamma_i,0)$ is expressed as
\begin{align}
\ddot x_1 - mgl_1 x_1 = 0, \label{eq:x1_inv_red}\\
\ddot x_2 - mgl_2 x_2 = 0, \label{eq:x2_inv_red}\\
\dot x_3 = 0. \label{eq:x3_inv_red}
\end{align}
Note that the inverted equilibrium of the Lagrange--Poincar\'e reduced equations has two negative eigenvalues, two positive eigenvalues and a zero eigenvalue. Thus, the inverted equilibrium $(\Gamma_i,0)$ is unstable and locally there exists a two dimensional stable manifold, a two dimensional unstable manifold and a one dimensional center manifold.

\begin{prop}\label{prop:inv_reduc}  
The inverted equilibrium $(\Gamma_i,0) = \left(-\frac{\rho}{\|\rho\|},0\right)$, of the Lagrange--Poincar\'e reduced dynamics of the 3D pendulum described by equations (\ref{eq:Reduc_Jw_dot}) and (\ref{eq:Gamma_dot}) is unstable.
\end{prop}

\begin{remark}{\rm Note that combining Proposition \ref{prop:inv_reduc} with Proposition \ref{prop:identif} immediately confirms the result that the inverted equilibrium manifold $\mathtt{I}$ of the full equations for the 3D pendulum given by (\ref{eq:Jw_dot})--(\ref{eq:R_dot}) is unstable.}
\end{remark}

We have analyzed the local stability properties of the hanging equilibrium and of the inverted equilibrium. We have not analyzed local stability properties of any other equilibrium solutions, but this analysis can easily be carried out using the methods that have been introduced.

\section{Lagrange--Routh Reduced 3D Pendulum Dynamics on $TS^2$}\label{sec:TS}
In the previous sections we studied the full and the Lagrange--Poincar\'e reduced equations of motion of the 3D pendulum. These involved the study of the dynamics of the 3D pendulum on $TSO(3)$ and on $TSO(3)/S^1$, respectively, using $(R,\omega)$ and $(\Gamma,\omega)$ to express the equations of motion. In this section, we present Lagrange--Routh reduction of the 3D pendulum, and we study the equations of motion that describe the evolution of $(\Gamma, \dot\Gamma) \in TS^2$.

\subsection{Lagrange--Routh Reduction of the 3D pendulum}

The key feature of Lagrange--Routh reduction is reducing the configuration
space into the quotient space induced by the symmetry action. The resulting
equations of motion on the reduced space are described in terms of the Euler-Lagrange
equation, but not with respect to the Lagrangian itself but with respect to the
Routhian~\cite{marsden_reduc, marsden_geo_top, holmes_marsden}.

The 3D pendulum has a $S^1$ symmetry given by a rotation about the vertical
axis. The symmetry action $\Phi_\theta: SO(3) \rightarrow SO(3)$ is given by
\begin{align*}
\Phi_{\theta}(R)=\exp(\theta\widehat{e}_3)R,
\end{align*}
for $\theta\in S^1$ and $R\in SO(3)$. It can be shown that the Lagrangian of the
3D pendulum is invariant under this symmetry action. Thus, the configuration
space is reduced to the shape manifold $SO(3) /S^1 \cong S^2$, and the dynamics
of the 3D pendulum is described in the tangent bundle $TS^2$. This reduction
procedure is interesting and challenging, since the projection $\Pi:SO(3)
\rightarrow S^2$ given by $\Pi(R)=R^T e_3$ together with the symmetry action has a nontrivial principal bundle structure. In other words, the angle of the rotation about
the vertical axis is not a global cyclic variable.

Here we present expressions for the Routhian and the reduced equations of
motion. The detailed description and development can be found
in the Appendix.

\begin{prop}[\cite{marsden_reduc}]
We identify the Lie algebra of $S^1$ with $\Re$. For $(R,\omega)\in T_R SO(3)$, the momentum map $\J:T SO(3) \rightarrow\Re^*$, the
locked inertia tensor $\mathbb{I}(R):\Re\rightarrow\Re^*$, and the mechanical
connection $\mathcal{A}:T SO(3) \rightarrow \Re$ for the 3D pendulum are given as
follows
\begin{align}
    \J(R,\hat\omega)&=e_3^T R J\omega,\\
    \mathbb{I}(R) &= e_3^T RJR^T e_3,\\
    \mathcal{A}(R,\hat\omega) & = \frac{e_3^T RJ\omega}{e_3^T RJR^T
    e_3}.\label{eqn:mechcon}
\end{align}
\end{prop}
\noindent The value of the momentum map $\mu=\J(R,\hat\omega)$ corresponds to
the vertical component of the angular momentum. Noether's theorem states that
the symmetry of the Lagrangian implies conservation of the corresponding momentum map. This is an alternative method of showing the invariant properties of the 3D pendulum, as opposed to the direct computation used in Section \ref{ssec:TSOinvar}.

Based on the above expressions, Lagrange--Routh reduction is carried out to
obtain the following result.
\begin{prop}\label{prop:routhred3Dpend}
For a given value of the momentum map $\mu$, the Routhian of the 3D pendulum is
given by
\begin{align}
R^{\mu}(\Gamma,\dot\Gamma) & = \frac{1}{2} (\dot\Gamma\times\Gamma)\cdot
J(\dot\Gamma\times\Gamma)-\frac{1}{2} (b^2+\nu^2) (\Gamma\cdot J\Gamma)
+mg\Gamma\cdot\rho,\label{eqn:Routhian}
\end{align}
where $b=\frac{ J\Gamma\cdot(\dot\Gamma\times\Gamma)}{\Gamma\cdot J\Gamma}$,
$\nu=\frac{\mu}{\Gamma\cdot J\Gamma}$, and the magnetic two form can be written
as
\begin{align}
\beta_\mu(\Gamma\times\eta,\Gamma\times\zeta) & = -\frac{\mu}{(\Gamma\cdot  J
\Gamma)^2} \bracket{ -(\Gamma\cdot  J \Gamma)\tr{J}
+2\norm{J\Gamma}^2}\Gamma\cdot(\eta\times\zeta).\label{eqn:betamu}
\end{align}
The Routhian satisfies the Euler-Lagrange equation, with the magnetic term,
given by
\begin{align}
\delta \int_{0}^T R^\mu (\Gamma,\dot\Gamma) dt = \int_{0}^T
\mathbf{i}_{\dot\Gamma} \beta_\mu(\delta\Gamma) dt.\label{eqn:EL}
\end{align}
This yields the reduced equation of motion on $TS^2$:
\begin{align}
\ddot\Gamma=-\|\dot\Gamma\|^2\Gamma+\Gamma\times\Sigma,\label{eqn:Routh_Reduced}
\end{align}
where
\begin{align}
\Sigma &= b\dot\Gamma +
J^{-1}\bracket{(J(\dot\Gamma\times\Gamma)-bJ\Gamma)\times ((\dot\Gamma\times
\Gamma)-b\Gamma)
+\nu^2 J\Gamma\times \Gamma -mg\Gamma\times\rho-c\dot\Gamma},\label{eq:Sigma}\\
c&=\nu\braces{\tr{J}-2\frac{\norm{J\Gamma}^2}{\Gamma\cdot J\Gamma}},\quad
b=\frac{ J\Gamma\cdot(\dot\Gamma\times\Gamma)}{\Gamma\cdot J\Gamma},\quad
\nu=\frac{\mu}{\Gamma\cdot J\Gamma}.\label{eqn:consts}
\end{align}
\end{prop}
\begin{proof} See the Appendix.\hfill\end{proof}

\subsection{Lagrange--Routh Reconstruction of the 3D pendulum}

For a given value of the momentum map $\mu$, let $\Gamma(t)\in S^2$ be a curve
in the reduced space $S^2$ satisfying the Euler-Lagrange equation for the
reduced Routhian $R^\mu$ given by \refeqn{Routh_Reduced}. The reconstruction
procedure is to find the curve $R(t)\in SO(3)$ in the configuration manifold that
satisfies $\Pi(R(t))=\Gamma(t)$ and $\J(R(t), R(t)^T\dot R(t))=\mu$.

This is achieved in two steps. First, we choose any curve $R_\hor(t)\in SO(3)$ in
the configuration manifold such that its projection is equal to the reduced
curve, i.e. $\Pi(R_\hor(t))=\Gamma(t)$. Now, the curve $R(t)$ can be written as
$R(t)=\Phi_{\theta(t)} (R_\hor (t))$ for some $\theta(t)\in S^1$. We find a
differential equation for $\theta(t)$ so that the value of the momentum map for
the reconstructed curve is conserved.

\begin{prop}
Suppose that the integral curve of the Lagrange--Routh reduced equation
\refeqn{Routh_Reduced} is given by $(\Gamma(t),\dot\Gamma(t))\in TS^2$, and the
value of the momentum map $\mu$ is known. The following procedure reconstructs
the motion of the 3D pendulum to obtain $(R(t),\omega(t))\in T SO(3)$ such that
$\Pi(R(t))=\Gamma(t)$ and $\J(R(t), \omega(t))=\mu$.
\begin{enumerate}
    \item Horizontally lift $\Gamma(t)$ to obtain $R_\hor (t)$ by
    integrating the following equation with $R_\hor(0)=R(0)$.
    \begin{align}
        \dot R_\hor (t) = R_\hor(t)\hat\omega_\hor(t),\label{eqn:Rhor}
    \end{align}
    where
    \begin{align}
        \omega_\hor(t)=\dot\Gamma(t)\times\Gamma(t)-b(t)\Gamma(t).\label{eq:Hor_vel}
    \end{align}
    \item Determine $\theta_\dyn(t)\in S^1$ by the following
    equation.
    \begin{align}
        \theta_\dyn (t)=\int_0^t \frac{\mu}{\Gamma(s)\cdot
        J\Gamma(s)}\,ds.\label{eqn:thetadyn}
    \end{align}
    \item Reconstruct the curve in $T SO(3)$.
    \begin{align}
        R(t)&=\Phi_{\theta_\dyn (t)}(R_\hor (t)) = \exp
        [\theta_\dyn(t)\hat e_3] R_\hor(t),\label{eqn:Rreconst}\\
        \omega(t)&=\omega_\hor(t) + \nu(t)\Gamma(t).\label{eqn:w_rec}
    \end{align}
\end{enumerate}
\end{prop}
\begin{proof} See the Appendix.\hfill\end{proof}

This leads to the geometric phase formula that expresses the rotation angle
about the vertical axis along a closed integral curve of the reduced equation.
\begin{prop}
Let $\Gamma(t)$ be a closed curve in $S^2$, i.e. $\Gamma(0)=\Gamma(T)$ for some
$T$. The geometric phase $\theta_\geo(T)\in S^1$ of the 3D pendulum is defined by
the relationship $R(T)=\Phi_{\theta_\geo(T)}(R(0))$ when $\mu=0$. This can be
written as
\begin{align}\label{eqn:geopen}
    \theta_\geo(T)=\int_\mathcal{B}\frac{2\norm{J\Gamma(t)}^2-\tr{J}(\Gamma(t)\cdot
    J\Gamma(t))}{(\Gamma(t)\cdot J\Gamma(t))^2}\,dA,
\end{align}
where $\mathcal{B}$ is a surface in $S^2$ with $\Gamma(t)$ as boundary.
\end{prop}

\subsection{Invariants of the Lagrange--Routh Reduced Model}

In this section we find an invariant of the motion for the Lagrange--Routh reduced model
of the 3D pendulum, namely the total energy of the system. Note that the Lagrange--Routh
reduced equations of motion are derived by elimination of the conserved
vertical component of the body-fixed angular momentum. In later sections, we
make use of the constant energy surfaces to understand the dynamics of the 3D
pendulum.

\begin{prop}
The total energy
\begin{align}
    E=\frac{1}{2}(\dot\Gamma\times\Gamma+(\nu-b)\Gamma)^TJ(\dot\Gamma\times\Gamma+(\nu-b)\Gamma)
    -mg\rho^T \Gamma\label{eqn:E_TS}
\end{align}
is constant along motions of the Lagrange--Routh reduced equations for the 3D pendulum.
\end{prop}
\begin{proof}
Substituting the reconstruction equations for the angular velocity
(\ref{eq:Hor_vel}), \refeqn{w_rec} into the total energy expression
\refeqn{E_TSOS}, we obtain \refeqn{E_TS}. The time derivative of the total
energy is given by
\begin{align*}
    \dot E=(\dot\Gamma\times\Gamma+(\nu-b)\Gamma)^TJ(\ddot\Gamma\times\Gamma+(\dot\nu-\dot b)\Gamma+(\nu-b)\dot\Gamma)
    -mg\rho^T \dot \Gamma.
\end{align*}
Substituting the reduced equation of motion \refeqn{Routh_Reduced} into the
above equation and rearranging, we can show that $\dot E=0$. \hfill\end{proof}

\subsection{Equilibria of the Lagrange--Routh Reduced Model}
The Lagrange--Routh reduced model can be considered as the Lagrange-Poincar\'{e} reduced model where the angular velocity is projected onto $T_{\Gamma} S^2$. Thus, the equilibria structure of the Lagrange-Routh reduced model is equivalent to the Lagrange-Poincar\'{e} reduced model, but it is represented by conditions on the reduced attitude $\Gamma_{\rm e}$ and the value of the momentum map $\mu$ instead of $(\Gamma_{\rm e},\omega_{\rm e})$.

More explicitly, we study the equilibria structure of the Lagrange-Routh reduced model using (\ref{eqn:Routh_Reduced}), and we show it is equivalent to the families of equilibria presented in Proposition \ref{prop:equi_reduced}.

\begin{prop}\label{prop:equi_Routhred_stat}
Consider the Lagrange--Routh reduced model of the 3D pendulum given by
(\ref{eqn:Routh_Reduced}). The equilibria $(\Gamma_{\rm e},0)\in TS^2$ of the Lagrange-Routh model are given for $\mu\in\Re$ as follows.
\begin{enumerate}
\item The hanging equilibrium: $\parenth{\frac{\rho}{\norm{\rho}},0}\quad\mu=0,$
\item The inverted equilibrium: $\parenth{-\frac{\rho}{\norm{\rho}},0}\quad\mu=0,$
\item Two relative equilibria:
\begin{align}
\parenth{-\frac{J^{-1}\rho}{\norm{J^{-1}\rho}},0}\quad%
\mu=\pm\sqrt{\frac{mg}{\norm{J^{-1}\rho}^3}}\rho^T J^{-1}\rho,\label{eqn:RelEqu3r}
\end{align}
\item One dimensional relative equilibrium manifolds:
\begin{align}
    \parenth{-\mathrm{sgn}(\alpha)\frac{n_\alpha}{\norm{n_\alpha}},0}\quad%
    \mu=\pm\mathrm{sgn}(\alpha)\sqrt{\frac{mg}{\norm{n_\alpha}^3}}n_\alpha^T J n_\alpha,\label{eqn:RelEqu41r}
\end{align}
where $n_\alpha= (J-\frac{1}{\alpha}I_{3 \times 3})^{-1}\rho\in\Re^3$ corresponding to $\alpha \in \mathcal{L}_i$, $i \in \{1,2,3,4\}$ and
\[
\mathcal{L}_1 = (-\infty,0)\cup (\tfrac{1}{J_3},\infty),\ \mathcal{L}_2 = (0,\tfrac{1}{J_1}),\ \mathcal{L}_3 = (\tfrac{1}{J_1},\tfrac{1}{J_2}),\ \mathcal{L}_4 = (\tfrac{1}{J_2},\tfrac{1}{J_3}).
\]
Also, $\mathrm{sgn}(\cdot)$ denotes the sign function.

\end{enumerate}
The families of equilibria given in (\ref{eqn:RelEqu41r}) converge to the hanging equilibrium and the inverted equilibrium when $\alpha\rightarrow 0$, and they converge to the third equilibria given in (\ref{eqn:RelEqu3r}), when $\alpha\rightarrow\pm\infty$. If the vector from the pivot to the center of mass $\rho$, lies on a principal axis, i.e. $\rho\times e_i=0$ for some $i\in\{1,2,3\}$, then (\ref{eqn:RelEqu41r}) can be rewritten as $\braces{(e_i,0 ),(-e_i,0)}$ for any $\mu\in\Re$.
Furthermore, there exist additional equilibria under the following assumptions on the moment of inertia matrix $J$ and the vector from the pivot to the center of mass $\rho = [\rho_1\ \rho_2\ \rho_3]\T$.
\begin{enumerate}\setcounter{enumi}{4}
\item $J_1$, $J_2$ and $J_3$ are distinct and $\rho_i = 0$. Then there exist one-dimensional relative equilibrium manifolds:
\begin{align}
\parenth{-\dfrac{p_i}{\|p_i\|},0}\quad%
\mu=\pm\sqrt{\dfrac{mg}{\|p_i\|^3}}p_i^T Jp_i,
\end{align}
for $i \in \{1,2,3\}$, where %
$p_1 = \big(\gamma,\frac{\rho_2}{J_2-J_1},\frac{\rho_3}{J_3-J_1}\big)$, %
$p_2 = \big(\frac{\rho_1}{J_1-J_2},\gamma,\frac{\rho_3}{J_3-J_2}\big)$, %
$p_3 = \big(\frac{\rho_1}{J_1-J_3},\frac{\rho_2}{J_2-J_3},\gamma\big)$ and $\gamma \in \mathbb{R}$.
\item $J_1 = J_2 \not= J_3$.
    \begin{enumerate}
    \item If $\rho_1 = \rho_2 = 0$. Then there exist two-dimensional relative equilibrium manifolds:
    \[
    \left(-\dfrac{q}{\|q\|},0\right)\quad%
    \mu=\pm\sqrt{\dfrac{mg}{\|q\|^3}}q^T Jq,
    \]
    where $q = \big(\gamma,\delta,\frac{\rho_3}{J_3-J_1}\big)$ and $\gamma,\delta \in \mathbb{R}$.
    \item If $\rho_3 = 0$. Then there exist one-dimensional relative equilibrium manifolds:
    \[
    \left(-\dfrac{r}{\|r\|},0\right)\quad%
    \mu=\pm\sqrt{\dfrac{mg}{\|r\|^3}}r^T Jr,
    \]
    where $r = \big(\frac{\rho_1}{J_1-J_3},\frac{\rho_2}{J_1-J_3},\gamma\big)$ and $\gamma \in \mathbb{R}$.
    \end{enumerate}
\end{enumerate}
\end{prop}

\begin{proof}
Substituting $\dot\Gamma_{\rm e}=0$ into \refeqn{Routh_Reduced}--\refeqn{consts}, we obtain a condition for an equilibrium $\Gamma_{\rm e}$ for $\mu$ as
\begin{align*}
    \Gamma_{\rm e} \times J^{-1}\bracket{\nu^2J\Gamma_{\rm e} \times \Gamma_{\rm e}-mg \Gamma_{\rm e} \times \rho }=0.
\end{align*}
This is equivalent to
\begin{align}
    \bracket{\nu^2J\Gamma_{\rm e} \times \Gamma_{\rm e}-mg \Gamma_{\rm e} \times \rho } = k_2 J\Gamma_{\rm e}\label{eqn:LRequcond0}
\end{align}
for a constant $k_2\in\Re$. Taking the dot product of this and $\Gamma_{\rm e}$ implies that $0=k_2 \Gamma_{\rm e}^T  J\Gamma_{\rm e}$. Since $\Gamma_{\rm e}^T  J\Gamma_{\rm e}>0$ as the moment of inertia matrix $J$ is positive definite and $\Gamma_{\rm e}\in S^2$, it follows that $k_2=0$. Thus, \refeqn{LRequcond0} is equivalent to
\begin{align}
    \nu^2 J\Gamma_{\rm e} + mg \rho = k_1\Gamma_{\rm e}\label{eqn:LRequcond}
\end{align}
for a constant $k_1\in\Re$. Note that this is equivalent to the equilibrium condition for the Lagrange-Poincar\'{e} reduced model given by \refeqn{cond0}: for any solution $(\Gamma_{\rm e},k,k_1)$ of \refeqn{cond0}, we can choose $\mu$ such that
$k^2=\nu^2=\frac{\mu^2}{(\Gamma_{\rm e}^T  J\Gamma_{\rm e})^2}$, which gives a solution of \refeqn{LRequcond}, and vice versa. Thus, the equilibria structure of the Lagrange-Routh reduced model is equivalent to the equilibria of the Lagrange-Poincar\'{e} reduced model. For an equilibrium $(\Gamma_{\rm e},\omega_{\rm e})$ of the Lagrange-Poincar\'{e} reduced model, the value of the momentum map at the corresponding equilibrium of the Lagrange-Routh model is given by
\begin{align}
    \mu = k (\Gamma_{\rm e}^T J\Gamma_{\rm e}) = \omega_{\rm e}^T \Gamma_{\rm e} (\Gamma_{\rm e}^T J\Gamma_{\rm e})
\end{align}
Substituting this into the equilibria presented in Proposition \ref{prop:equi_reduced}, we obtain the equilibria of the Lagrange-Routh reduced model.\hfill\end{proof}

\subsection{Local Analysis of the Lagrange--Routh Reduced Model on $TS^2$}

We showed that in case $\mu = 0$, the Routh reduced model of the 3D pendulum has two isolated equilibria, namely the hanging equilibrium and the inverted equilibrium. These equilibria correspond to the disjoint equilibrium manifolds of the full equations of the 3D pendulum.

We next focused on these isolated equilibria of the Routh reduced equations. Using Proposition \ref{prop:equi_Routhred_stat}, the stability properties of the equilibrium manifolds of the 3D pendulum can be deduced by studying the Routh reduced equilibria for the case $\mu =0$. Compared to the Lagrange--Poincar\'e reduced model, the Routh reduction procedure results in a set of complicated equations that are a challenge to analyze.

We now present local analyses of the Routh reduced model of the 3D pendulum near the hanging equilibrium and near the inverted equilibrium.

Consider the equations (\ref{eq:x1_hang})--(\ref{eq:x3_hang}) representing the linearization of the full equations of motion of the 3D pendulum at the hanging equilibrium. It was shown before that the Lagrange--Poincar\'e reduced equations of motion can be written in terms of $(x_1,x_2,\dot x_1, \dot x_2, \dot x_3)$. As shown in Proposition \ref{prop:Reduc_3D_Pend_linear_hng}, this result follows from the fact that any perturbation in $\Gamma \in S^2$ at $\Gamma_h$ can be expressed in terms of $(x_1,x_2) \in \mathbb{R}^2$. In a similar fashion, one obtains the following result.

\begin{prop}\label{prop:Routh_Reduc_3D_Pend_linear_hng}
The linearization of the Routh reduced attitude dynamics of the 3D pendulum, at the equilibrium $(\Gamma_h,0) = (R_{\rm e}\T e_3,0)$, described by equation (\ref{eqn:Routh_Reduced}) can be expressed using $(x_1,x_2,\dot x_1, \dot x_2) \in \mathbb{R}^4$ according to (\ref{eq:x1_hang}) and (\ref{eq:x2_hang}).
\end{prop}
\begin{proof} In Proposition \ref{prop:Reduc_3D_Pend_linear_hng}, it was shown that the perturbations in $\Gamma$ at $\Gamma_h$ can be described in terms of $(x_1,x_2) \in \mathbb{R}^2$. The proof then simply follows by noting that the equations of motion of the Routh reduced 3D pendulum is described in terms of $(\Gamma,\dot \Gamma) \in TS^2$. Thus the linearization of the  Routh reduced 3D pendulum model can be described using using $(x_1,x_2,\dot x_1, \dot x_2) \in \mathbb{R}^4$ according to (\ref{eq:x1_hang}) and (\ref{eq:x2_hang}).\hfill\end{proof}

The linearization of the Routh reduced attitude dynamics of the 3D pendulum, about the equilibrium $(0,\Gamma_h)$ is obtained from the linearized model of the full attitude dynamics by neglecting the dynamics corresponding to $x_3$. Summarizing the above, the linearization of (\ref{eqn:Routh_Reduced}) about the hanging equilibrium $(0,\Gamma_h)$ is expressed as
\begin{align}
\ddot x_1 + mgl_1 x_1 = 0, \label{eq:x1_hang_routh_red}\\
\ddot x_2 + mgl_2 x_2 = 0. \label{eq:x2_hang_routh_red}
\end{align}
It is clear that due to the presence of imaginary eigenvalues, stability of the hanging equilibrium $(\Gamma_h,0)$ cannot be concluded. Therefore we next consider Lyapunov analysis.

\begin{prop}\label{prop:hang_Routh_reduc}  The hanging equilibrium $(\Gamma_h,0) = \left(\frac{\rho}{\|\rho\|},0\right)$ of the reduced dynamics of the 3D pendulum described by (\ref{eqn:Routh_Reduced}) is stable in the sense of Lyapunov.
\end{prop}
\begin{proof} Consider the candidate Lyapunov function
\begin{equation}\label{energy_routh_reduc}
V(\Gamma, \dot\Gamma) =\frac{1}{2}(\dot\Gamma\times\Gamma+(\nu-b)\Gamma)\T J(\dot\Gamma\times\Gamma+(\nu-b)\Gamma) + mg(\|\rho\|-\rho\T \Gamma).
\end{equation}
Note that $V(\Gamma_h,0)=0$ and $V(\Gamma,\dot\Gamma) > 0$ elsewhere. Furthermore, the derivative along a solution of (\ref{eq:Reduc_Jw_dot}) and (\ref{eq:Gamma_dot}) is given by
\[
\begin{aligned}
  \dot{V}(\Gamma, \dot\Gamma) &= (\dot\Gamma\times\Gamma+(\nu-b)\Gamma)\T J(\ddot\Gamma\times\Gamma+(\dot\nu-\dot b)\Gamma+(\nu-b)\dot\Gamma) -mg\rho\T \dot \Gamma.
\end{aligned}
\]
Substituting the reduced equation of motion (\ref{eqn:Routh_Reduced}) into the above equation and rearranging, we can show that $\dot{V}(\Gamma, \dot\Gamma) =0$. Thus, the hanging equilibrium of (\ref{eqn:Routh_Reduced}) is Lyapunov stable.\hfill\end{proof}

\begin{remark}{\rm Note that combining Proposition \ref{prop:hang_Routh_reduc} with Proposition \ref{prop:equi_Routhred_stat} immediately yields the result in Proposition \ref{thm:hang}.}
\end{remark}

We next study the local properties of the Routh reduced equations of the 3D pendulum near the inverted equilibrium $(\Gamma_i,0)$. Consider the linearization of (\ref{eqn:Routh_Reduced}) at an equilibrium $(\Gamma_i,0)= (R_{\rm e}\T e_3,0)$, where $(R_{\rm e},0)$ is an equilibrium of the inverted equilibrium manifold $\mathtt{I}$. A result similar to Proposition \ref{prop:Routh_Reduc_3D_Pend_linear_hng} follows.

\begin{prop}\label{prop:Rout_Reduc_3D_Pend_linear_inv}
The linearization of the reduced attitude dynamics of the 3D pendulum, at the equilibrium $(\Gamma_i,0) = (R_{\rm e}\T e_3,0)$, described by (\ref{eqn:Routh_Reduced}) can be expressed using $(x_1,x_2,\dot x_1, \dot x_2) \in \mathbb{R}^4$ according to (\ref{eq:x1_inv}) and (\ref{eq:x2_inv}).
\end{prop}

Summarizing the above, the linearization of (\ref{eqn:Routh_Reduced}) at the inverted equilibrium $(\Gamma_i,0)$ is expressed as
\begin{align}
\ddot x_1 - mgl_1 x_1 = 0, \label{eq:x1_inv_Routh_red}\\
\ddot x_2 - mgl_2 x_2 = 0. \label{eq:x2_inv_Routh_red}
\end{align}
Note that the linearization of (\ref{eqn:Routh_Reduced}) at the inverted equilibrium has two negative eigenvalues and two positive eigenvalues. Thus, the inverted equilibrium $(\Gamma_i,0)$ of the Routh reduced model is unstable and locally there exists a two dimensional stable manifold and a two dimensional unstable manifold.

\begin{prop}\label{prop:inv_Routh_reduc} 
The inverted equilibrium $(\Gamma_i,0) = \left(-\frac{\rho}{\|\rho\|},0\right)$ of the Routh reduced dynamics of the 3D pendulum described by (\ref{eqn:Routh_Reduced}) is unstable.
\end{prop}

\begin{remark}{\rm Note that combining Proposition \ref{prop:inv_Routh_reduc} with Proposition \ref{prop:equi_Routhred_stat} immediately yields the result that the inverted equilibrium manifold $\mathtt{I}$ of the 3D pendulum given by (\ref{eq:Jw_dot})--(\ref{eq:R_dot}) is unstable.}
\end{remark}

\subsection{Poincar\'{e} Map on the Lagrange--Routh Reduced Model}

A Poincar\'{e} map describes the evolution of intersection points of a trajectory with a transversal hypersurface of codimension one. Typically, one chooses a hyperplane, and considers a trajectory with initial conditions on the hyperplane. The points at which this trajectory returns to the hyperplane are then observed, which provides insight into the stability of periodic orbits or the global characteristics of the dynamics.

The Lagrange--Routh reduced equations for the 3D pendulum on $TS^2$ are a particularly suitable choice for analysis using a Poincar\'{e} map, since it has dimension 4.  Since the total energy given by 
(\ref{eqn:E_TS}) is conserved, choosing a Poincar\'e section on $TS^2$ and restricting to an energy isosurface induces a Poincar\'e map on a 2-dimensional submanifold of $TS^2$. We define a Poincar\'{e} section on $TS^2$ for the Lagrange-Routh dynamics of the 3D pendulum given by (\ref{eqn:Routh_Reduced}) as follows.
\begin{align*}
    \mathcal{P}=\braces{(\Gamma,\dot\Gamma)\in TS^2\,\big|\,e_3^T\dot\Gamma=0,\;
    e_3^T(\Gamma\times\dot\Gamma) > 0, \text{ and $E(\Gamma,\dot\Gamma) =$ constant }}.
\end{align*}
Suppose $\Gamma\in\mathcal{P}$ is given. The tangent space $T_\Gamma S^2$ is a
plane that is tangential to $S^2$ and perpendicular to $\Gamma$. The first
condition of the Poincar\'{e} section, $e_3^T\dot\Gamma=0$, determines a line
in which the tangent vector $\dot\Gamma\in T_\Gamma S^2$ should lie, and the
constraint of the total energy conservation fixes the magnitude of the tangent
vector in that line. Thus, the tangent vector is uniquely determined up to sign
change. The second condition of the Poincar\'{e} section resolves this
ambiguity. It also excludes two reduced attitudes $\Gamma=\pm e_3$ for which
the first condition is trivial; $e_3^T\dot\Gamma=0$ for any $\dot\Gamma\in
T_{e_3}S^2\bigcup T_{-e_3}S^2$. Thus, $\mathcal{P}$ can be equivalently identified as
\begin{align*}
    \mathcal{P}=\braces{\Gamma\in S^2\,\big|\,e_3^T\dot\Gamma=0,\;
    e_3^T(\Gamma\times\dot\Gamma) > 0, \text{ and $E(\Gamma,\dot\Gamma) =$ constant }},
\end{align*}
where $(\Gamma,\dot\Gamma)$ satisfies (\ref{eqn:Routh_Reduced}).

This Poincar\'{e} section in $TS^2$ is well-defined in the sense that for each
element, the corresponding tangent vector is uniquely determined. The
attitude and the angular velocity in $T SO(3)$ can be obtained by using the
reconstruction procedure for the given value of the momentum map $\mu$. 

\begin{figure}
\centerline{
    \subfigure[$E=-2.65$]{
    \includegraphics[width=0.2\textwidth]{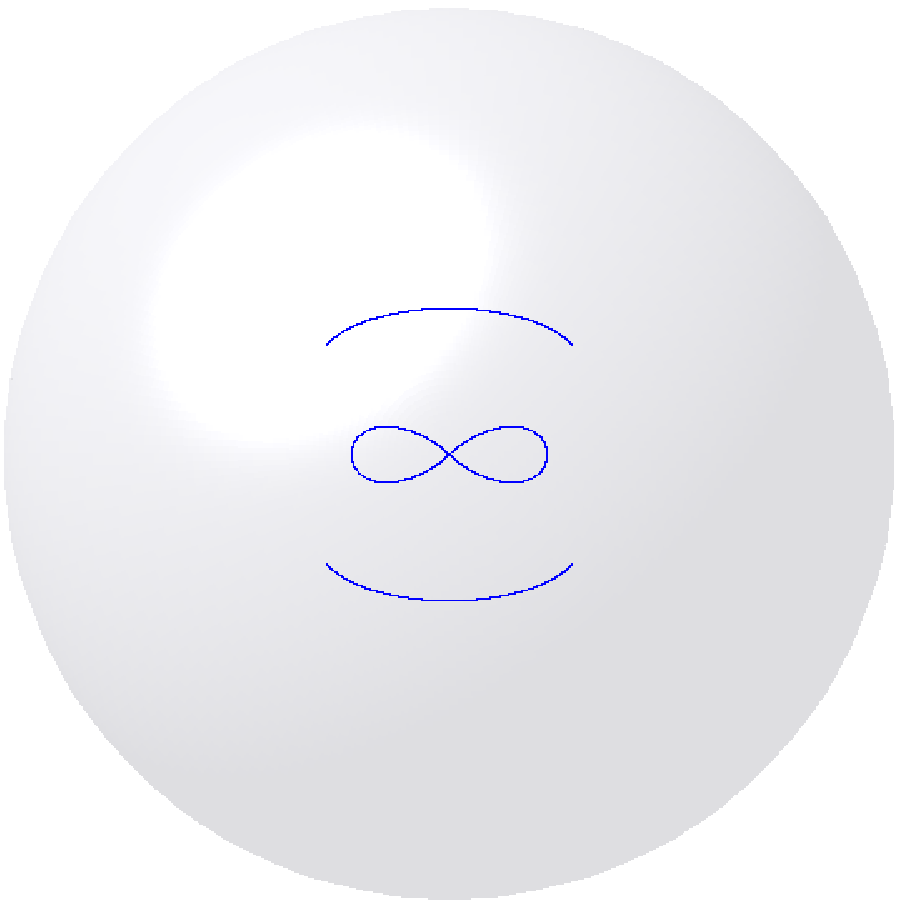}\label{fig:poin1}}
    \subfigure[$E=0$]{
    \includegraphics[width=0.2\textwidth]{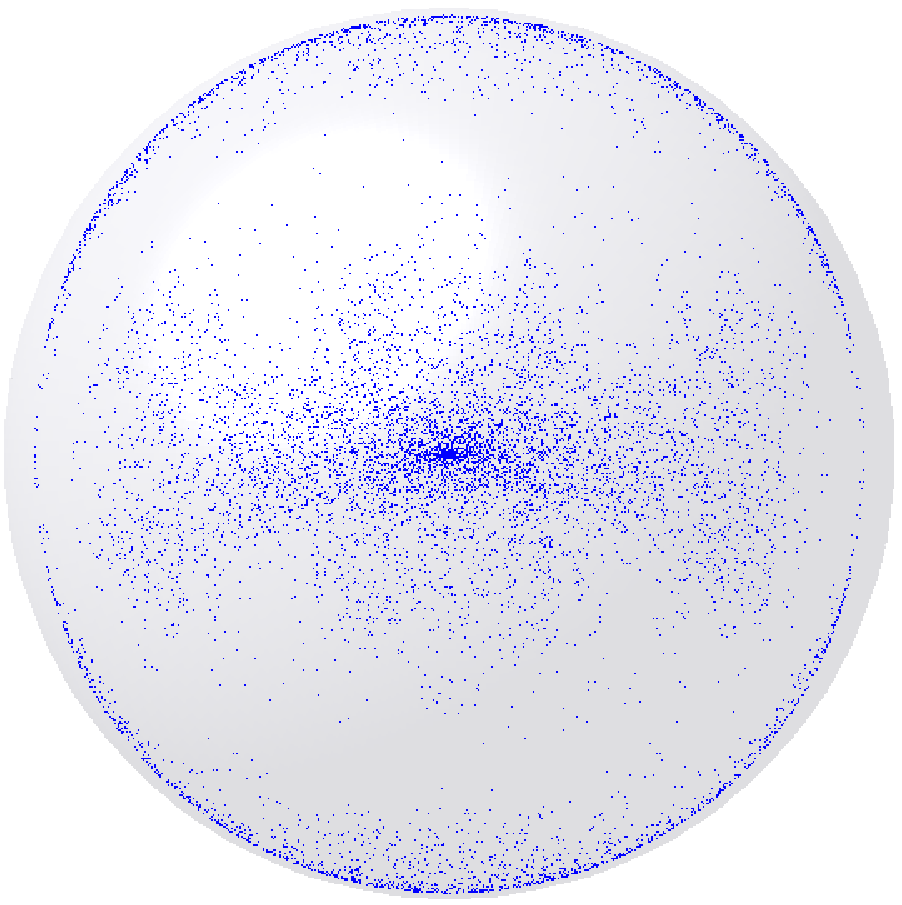}\label{fig:poin2}}
    \subfigure[$E=2.03$]{
    \includegraphics[width=0.2\textwidth]{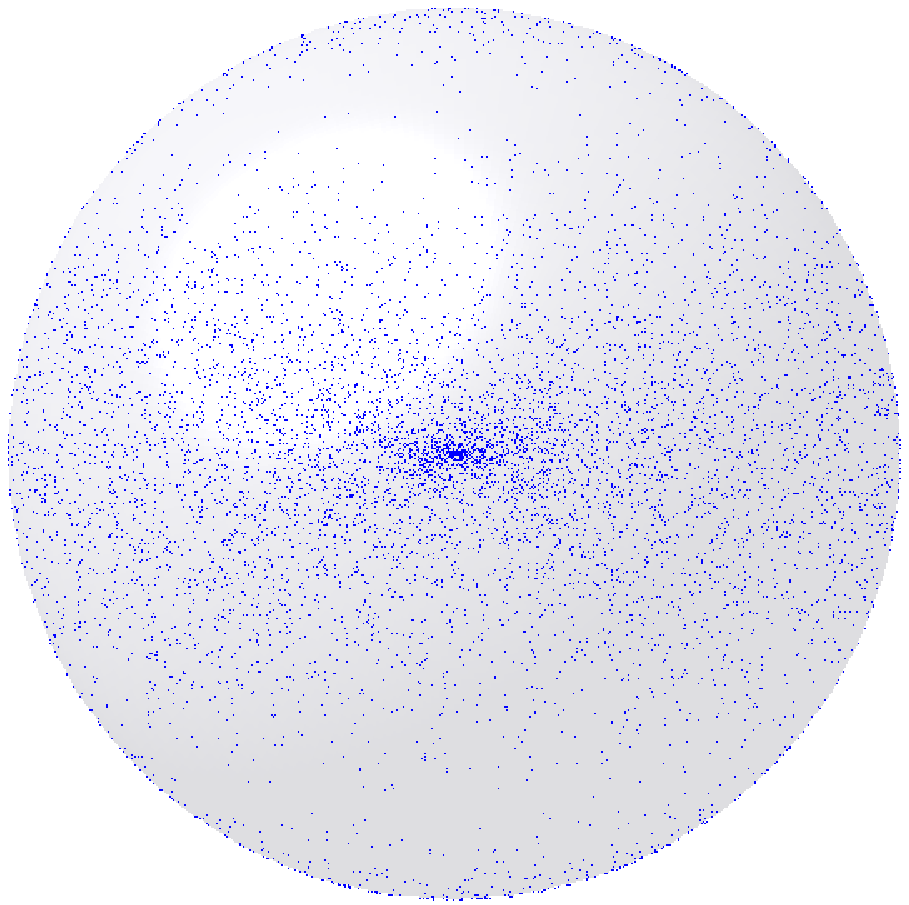}\label{fig:poin3}}
    \subfigure[$E=8.83$]{
    \includegraphics[width=0.2\textwidth]{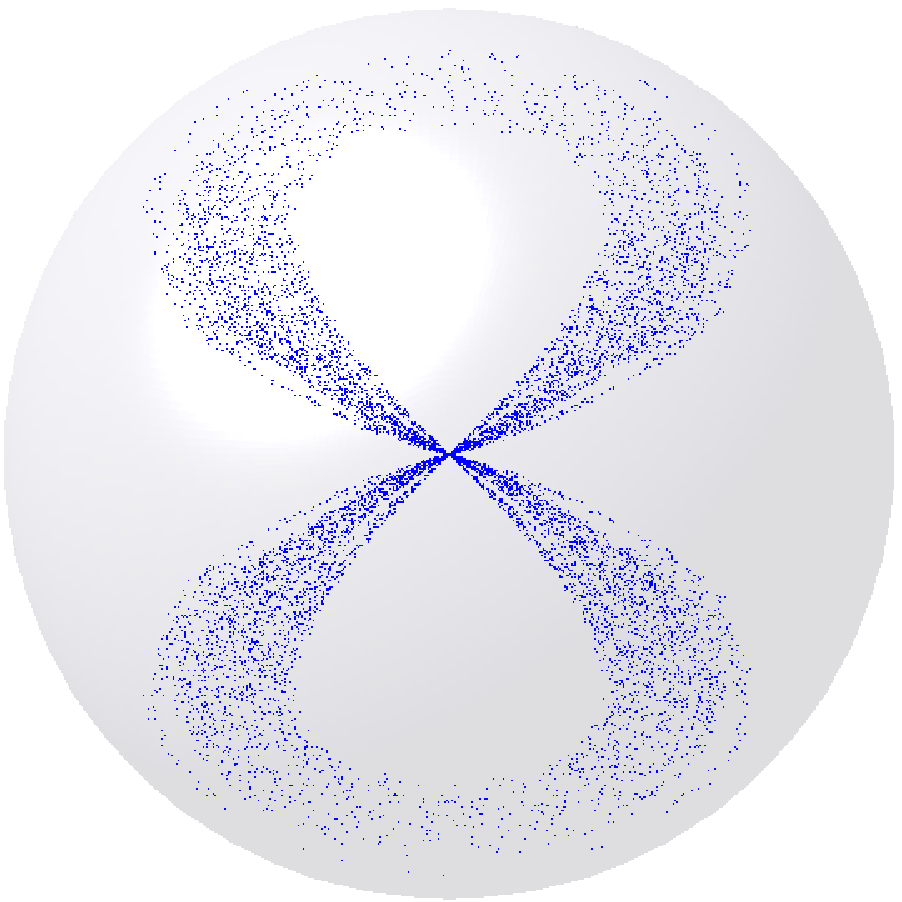}\label{fig:poin4}}
    \subfigure[$E=11.95$]{
    \includegraphics[width=0.2\textwidth]{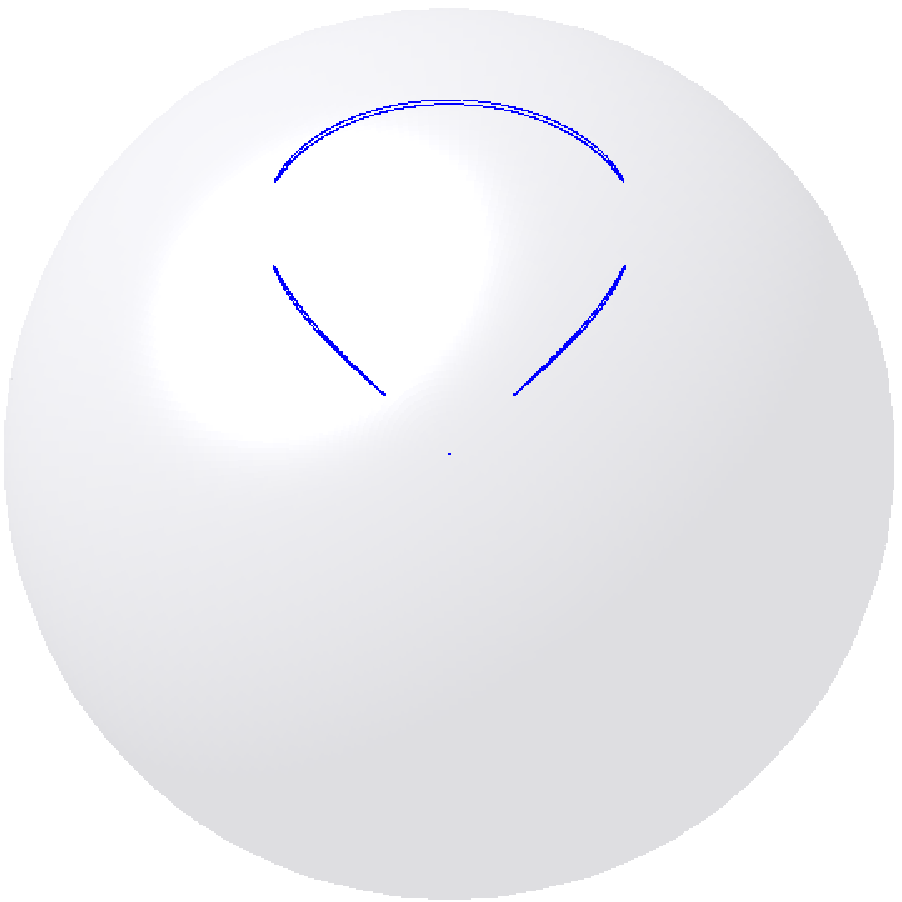}\label{fig:poin5}}
    }
\caption{Poincar\'{e} maps for 3D pendulum with varying total
energy}\label{fig:poin}
\end{figure}

\reffig{poin} shows particular examples for this Poincar\'{e} map. The pendulum
body is chosen as an elliptic cylinder with properties of $m=1\,\mathrm{kg}$,
$J=\mathrm{diag}[0.13,0.28,0.17]\,\mathrm{kgm^2}$,
$\rho=[0,0,0.3]\,\mathrm{m}$. The initial condition are given by $R_0=I_{3\times
3}$ and $\omega_0=c[1,1,1]\,\mathrm{rad/s}$, where the constant $c$ is varied
to give different total energy levels. The Lie group variational integrator
introduced in \cite{CCA05} is used to compute the Poincar\'{e} maps numerically.

It is interesting to see the transition of the Poincar\'{e} maps with varying
total energy levels. The attitude dynamics of the 3D pendulum is periodic in
\reffig{poin1}, but it exhibits chaotic behavior with increased energy level
in \reffig{poin2} and \ref{fig:poin3}. If the total energy is increased
further, the attitude dynamics becomes periodic again in \reffig{poin5}. This
demonstrates the highly-nonlinear, and perhaps chaotic, characteristics of the
3D pendulum dynamics.

\section{Conclusions}

The 3D pendulum exhibits rich dynamics with nontrivial geometric structure; these dynamics are much richer and more complex than the dynamics of a 1D planar pendulum or a 2D spherical pendulum.   This paper has demonstrated that the methods of geometric mechanics and the methods of nonlinear dynamics can be meshed to obtain insight into the complex dynamics of the 3D pendulum.

We have introduced three different models of the 3D pendulum, including the full model defined on $TSO(3)$, the Lagrange--Poincar\'e reduced model on $TSO(3)/S^1$ obtained by identifying configurations in the same group orbit, and the Lagrange--Routh reduced model on $TS^2$ where one additionally utilizes the fact that the dynamics evolves on a momentum level set.  Relationships between the various representations are discussed in the context of conservation properties, equilibria and their stability properties, and invariant manifolds.

In addition, we illustrate that the use of the Lagrange--Routh reduced equations of motion, together with the energy conservation properties, allow the construction of a Poincar\'e map that can be readily visualized, thereby providing a graphical tool for obtaining insight into the rich nonlinear dynamical properties of the 3D pendulum.

\Appendix

\section*{Appendix}

In this appendix, we summarize Lagrange-Routh reduction and reconstruction procedures for the 3D pendulum.

\subsection{Reduction}

A description of Lagrange-Routh reduction can be found in~\cite{marsden_reduc} including expressions for the mechanical connection and the Routhian of the 3D pendulum given by \refeqn{mechcon} and \refeqn{Routhian}, respectively. Here we derive the reduced equation of motion \refeqn{Routh_Reduced} using the Euler-Lagrange equation for the given Routhian \refeqn{Routhian}.

\paragraph*{Variation of Routhian}
The Routhian satisfies the variational Euler-Lagrange
equation with the magnetic term given by \refeqn{EL}. We use constrained variations of $\Gamma\in S^2$:
\begin{align}
\delta\Gamma & = \Gamma \times \eta,\label{eqn:delGam}\\
\delta\dot\Gamma & = \dot\Gamma \times \eta + \Gamma \times \dot
\eta.\label{eqn:delGamd}
\end{align}
Here we assume that $\eta\cdot\Gamma=0$, since the component of $\eta$ parallel
to $\Gamma$ has no effect on $\delta\Gamma$. These expressions are essential for developing the reduced equation of motion.

Using \refeqn{delGam}, \refeqn{delGamd}, and the properties $\Gamma\cdot\dot\Gamma=0$, $\Gamma\cdot\eta=0$, the variation of the Routhian is given by
\begin{align}
    \delta R^\mu = \dot\eta\cdot J(\dot\Gamma\times\Gamma-b\Gamma)
    -\eta\cdot\Gamma\times\bracket{-\dot\Gamma\times
    J(\dot\Gamma\times\Gamma)+(b^2+\nu^2)J\Gamma-bJ(\dot\Gamma\times\Gamma)+b(\dot\Gamma\times J\Gamma)+mg\rho}.\label{eqn:delR}
\end{align}

\paragraph*{Magnetic 2-form}
From the given mechanical connection $\mathcal{A}$ and a value of the momentum map $\mu\in\Re^*$, define a 1-form $\mathcal{A}_\mu$ on $T SO(3)$ by
\begin{align*}
\mathcal{A}_\mu(R)\cdot (R,\hat\omega)=\pair{\mu,\mathcal{A}(R,\hat\omega)}
=\mu \frac{e_3^T RJ\omega}{e_3^T R J R^T e_3}
\end{align*}
The magnetic 2-form $\beta_\mu$ in \refeqn{betamu} is the exterior derivative of $\mathcal{A}_\mu$, which can be obtained by using the identity
$\mathbf{d}\mathcal{A}_{\mu} (X,Y) =
X[\mathcal{A}_{\mu}(Y)]-Y[\mathcal{A}_{\mu}(X)]-\mathcal{A}_{\mu}([X,Y])$
for $X=R\hat\eta,Y=R\hat\zeta\in T_R SO(3)$.
Suppose that $\dot\Gamma=\Gamma\times\omega$. Since
$\Gamma\cdot(\omega\times\eta)=\eta\cdot(\Gamma\times\omega)=\eta\cdot\dot\Gamma$,
the interior product of the magnetic 2-form is given by
\begin{align}
    \intp_{\dot\Gamma}\beta_\mu(\delta\Gamma)&=\beta_\mu(\Gamma\times\omega,\Gamma\times\eta)
     = \nu\braces{\tr{J}-2\frac{\norm{J\Gamma}^2}{\Gamma\cdot J\Gamma}}\dot\Gamma\cdot\eta,\label{eqn:ibeta}
\end{align}
where $\nu=\frac{\mu}{\Gamma\cdot J\Gamma}$.

\paragraph*{Euler-Lagrange equation with magnetic 2-form}
Substituting \refeqn{delR} and \refeqn{ibeta} into \refeqn{EL}, and integrating by parts, the Euler-Lagrange equation for the reduced Routhian \refeqn{Routhian} is written as
\begin{align}
-\int_0^T \eta\cdot\bracket{J(\ddot\Gamma\times \Gamma -b\dot\Gamma-\dot
b\Gamma)+\Gamma\times X+c\dot\Gamma}\, dt=0,\label{eqn:ELf}
\end{align}
where
\begin{align}
X=-\dot\Gamma\times J(\dot\Gamma\times
\Gamma)+(b^2+\nu^2)J\Gamma-bJ(\dot\Gamma\times\Gamma)+b(\dot\Gamma\times
J\Gamma)+mg\rho,\label{eqn:X}
\end{align}
and $c$ is given by \refeqn{consts}. Since \refeqn{ELf} is satisfied for all $\eta$ with $\Gamma\cdot\eta=0$, we obtain
\begin{align}
J(\ddot\Gamma\times \Gamma -b\dot\Gamma-\dot b\Gamma)+\Gamma\times
X+c\dot\Gamma=\lambda\Gamma,\label{eqn:eq0}
\end{align}
for $\lambda\in\Re$. This is the reduced equation of motion. However, this equation has an ambiguity since the value of $\lambda$ is unknown; this equation is implicit for $\ddot\Gamma$ since the term $\dot b$ is expressed in terms of $\ddot\Gamma$. The next step is to determine expressions for $\lambda$ and $\dot b$ using the definition of $b$ and some vector identities.

We first find an expression for $\lambda$ in terms of $\Gamma,\dot\Gamma$. Taking the dot product of \refeqn{eq0} with $\Gamma$, we obtain
\begin{align}
\Gamma\cdot J(\ddot\Gamma\times \Gamma -b\dot\Gamma-\dot
b\Gamma)=\lambda.\label{eqn:lam0}
\end{align}
From the definition of $b$, we can show the following identity: $\Gamma \cdot J(\dot\Gamma\times\Gamma -b \Gamma)=0$. Differentiating this with time, and substituting into \refeqn{lam0}, we find an expression for $\lambda$ in terms of $\Gamma,\dot\Gamma$ as
\begin{align}
\lambda=-\dot\Gamma \cdot J(\dot\Gamma\times\Gamma -b \Gamma).\label{eqn:lambda}
\end{align}
Substituting \refeqn{lambda} into \refeqn{eq0}, and taking the dot product of the result with $\Gamma$, we obtain an expression for $\dot b$ in terms of $\Gamma,\dot\Gamma$ as
\begin{align}
\dot b=\Gamma\cdot J^{-1}\braces{\Gamma\times X+c\dot\Gamma+(\dot\Gamma \cdot
J(\dot\Gamma\times\Gamma -b \Gamma))\Gamma}.\label{eqn:bdot}
\end{align}
Substituting \refeqn{bdot} into \refeqn{eq0}, and using the vector identity
$Y-(\Gamma\cdot Y)\Gamma=(\Gamma\cdot\Gamma)Y-(\Gamma\cdot
Y)\Gamma=-\Gamma\times(\Gamma\times Y)$ for any $Y\in\Re^3$, we obtain the
following form for the reduced equation of motion
\begin{align*}
\ddot\Gamma\times \Gamma -b\dot\Gamma-\Gamma\times\bracket{\Gamma\times J^{-1}
\braces{\Gamma\times X+c\dot\Gamma+(\dot\Gamma \cdot J(\dot\Gamma\times\Gamma
-b \Gamma))\Gamma}}=0.
\end{align*}

\paragraph*{Reduced equation of motion}
This equation has no ambiguity. Now, we simplify this equation. The above expression is equivalent to the following equation
\begin{align*}
\Gamma\times\bracket{\ddot\Gamma\times \Gamma
-b\dot\Gamma-\Gamma\times\bracket{\Gamma\times J^{-1}\braces{\Gamma\times
X+c\dot\Gamma+(\dot\Gamma \cdot J(\dot\Gamma\times\Gamma -b
\Gamma))\Gamma}}}=0.
\end{align*}
Since $\Gamma\cdot \ddot\Gamma=-\|\dot\Gamma\|^2$, the first term is given by
\begin{align*}
\Gamma\times(\ddot\Gamma\times\Gamma) & = (\Gamma\cdot\Gamma)\ddot\Gamma
-(\Gamma\cdot\ddot\Gamma)\Gamma=\ddot\Gamma+\|\dot\Gamma\|^2\Gamma.
\end{align*}
Using the property $\Gamma\times(\Gamma\times(\Gamma\times
Y))=-(\Gamma\cdot\Gamma)\Gamma\times Y=-\Gamma\times Y$ for $Y\in\Re^3$, the
third term of the above equation can be simplified. Substituting \refeqn{X} and rearranging, the reduced equation of motion for the 3D pendulum is given by
\begin{align}
\ddot\Gamma=-\|\dot\Gamma\|^2\Gamma+\Gamma\times\Sigma,
\end{align}
where $\Sigma = b\dot\Gamma+
J^{-1}\bracket{(J(\dot\Gamma\times\Gamma)-bJ\Gamma)\times ((\dot\Gamma\times
\Gamma)-b\Gamma)
+\nu^2 J\Gamma\times \Gamma -mg\Gamma\times\rho-c\dot\Gamma}$.

\subsection{Reconstruction}

For a given integral curve of the reduced equation $(\Gamma(t),\dot\Gamma(t))\in TS^2$, we find a curve $\tilde{R}(t)\in SO(3)$ that is projected into the reduced curve, i.e. $\Pi(\tilde{R}(t)=\Gamma(t)$. The reconstructed curve can be written as $R(t)=\Phi_{\theta(t)}( \tilde{R}(t))$ for some $\theta(t)\in S^1$. The conservation of the momentum map yields the following reconstruction equation~\cite{marsden_reduc}.
\begin{align*}
    \theta(t)^{-1}\dot\theta(t)=\mathbb{I}^{-1}(\tilde R(t))\mu - \mathcal{A}(\dot{\tilde R}(t)).
\end{align*}
The particular choice of $\tilde{R}(t)$, the horizontal lift given by \refeqn{Rhor}, simplifies the above equation, since the horizontal part of the tangent vector is annihilated by the mechanical connection. Further, since the group $S^1$ is abelian, the solution reduces to a quadrature such as \refeqn{thetadyn}. The reconstructed curve is given by \refeqn{Rreconst}.

\bibliographystyle{siam}
\bibliography{3Dpend}

\end{document}